\documentclass[a4paper,11pt,twoside,fleqn]{article}
\usepackage{amsmath,amssymb,amsthm,latexsym}
\newtheorem{theorem}{Theorem}[section]
\newtheorem{proposition}[theorem]{Proposition}
\newtheorem{corollary}[theorem]{Corollary}
\renewcommand{\Re}{\text{\rm Re}\,}
\renewcommand{\Im}{\text{\rm Im}\,}
\newcommand{\dist}{\text{\rm dist}}
\newcommand{\Id}{\text{\rm I}}

\begin{document}

\title{On solutions to the Ginzburg-Landau equations in higher dimensions}
\author{Simon Brendle}
\date{August 13, 2003}

\maketitle

\section{Introduction}

Let $M$ be a Riemannian manifold of dimension $n \geq 2$. Consider the semilinear elliptic 
equation 
\begin{equation} 
d^* d\phi = \frac{1}{2\varepsilon^2} \, (1 - |\phi|^2) \, \phi, 
\end{equation} 
where $\phi$ is a complex-valued function on $M$. This equation is the Euler-Lagrange equation 
for the functional 
\[E_\varepsilon(\phi) = \int_M \Big ( |d\phi|^2 + \frac{1}{4\varepsilon^2} \, (1 - |\phi|^2)^2 
\Big ).\] 
The equation (1) has been studied by many authors, including F. Bethuel, H. Brezis and 
F. H\'elein \cite{BBH}, F.-H. Lin \cite{Li1,Li2}, and R. Jerrard and H. M. Soner \cite{JS1,JS2}. 
The corresponding Schr\"odinger and wave equation were studied by J. Colliander and R. Jerrard 
\cite{CJ,Je} and by F.-H. Lin, J. Xin, and P. Zhang \cite{LX,LZ}. While these results are mainly 
devoted to the case $n = 2$, the higher-dimensional situation has been studied by F.-H. Lin and 
T. Rivi\`ere \cite{LR} and F. Bethuel, H. Brezis and G. Orlandi \cite{BBO}. \\

An important problem is to describe the behavior of the solutions as $\varepsilon \to 0$. 
Suppose that $\phi_j$ is a sequence of complex-valued functions on $M$ such that 
\[d^* d\phi_j = \frac{1}{2\varepsilon_j^2} \, (1 - |\phi_j|^2) \, \phi_j\] 
where $\varepsilon_j \to 0$. Then there exists a closed set $S$ of Hausdorff codimension $2$ and 
a harmonic map $\phi_\infty: M \setminus S \to S^1$ such that $\phi_j \to \phi_\infty$ on $M 
\setminus S$. In particular, if $M$ has dimension $2$, then the set $S$ is finite, and its 
cardinality is given by the degree of $\phi_j$. \\

In higher dimensions, it follows from results of F.-H. Lin and T. Rivi\`ere \cite{LR} and 
F. Bethuel, H. Brezis, and G. Orlandi \cite{BBO} that the vortex submanifold $S$ is stationary 
in the sense that its generalized mean curvature is equal to $0$. \\

In this paper, we study the converse problem. To this end, we consider a smooth minimal 
submanifold $S$ of codimension $2$. Our aim is to construct solutions of the Ginzburg-Landau 
equations 
\begin{equation} 
d^* F_A = \frac{1}{2\varepsilon^2} \, (\phi \, \overline{D_A \phi} - \overline{\phi} \, D_A 
\phi) 
\end{equation}
and 
\begin{equation} 
D_A^* D_A \phi = \frac{1}{2\varepsilon^2} \, (1 - |\phi|^2) \, \phi. 
\end{equation} 
Here, $A$ is a connection on a complex line bundle $L$ over $M$, and $\phi$ is a section of $L$. 
A pair $(A,\phi)$ satisfies (2), (3) if and only if $(A,\phi)$ is a critical point of the 
Ginzburg-Landau functional 
\[E_\varepsilon(A,\phi) = \int_M \Big ( \varepsilon^2 \, |F_A|^2 + |D_A \phi|^2 + \frac{1}{4
\varepsilon^2} \, (1 - |\phi|^2)^2 \Big ).\] 
In dimension $2$, Bogomol'nyi observed that the Ginzburg-Landau functional has a lower bound 
which depends only on the degree of the line bundle $L$. This is a consequence of the identity 
\[E_\varepsilon(A,\phi) = \int_{\mathbb{R}^2} \Big ( \varepsilon * \! (i F_A) - \frac{1}{2
\varepsilon} \, (1 - |\phi|^2) \Big )^2 + 2 \int_{\mathbb{R}^2} |\bar{\partial}_A \phi|^2 + 2\pi 
\, c_1(L).\] 
From this it follows that 
\[E_\varepsilon(A,\phi) \geq 2\pi \, c_1(L)\] 
with equality if and only if $(A,\phi)$ is a solution of the vortex equations 
\begin{equation} 
\varepsilon * \! (i F_A) = \frac{1}{2\varepsilon} \, (1 - |\phi|^2) 
\end{equation} 
and 
\begin{equation} 
\bar{\partial}_A \phi = 0. 
\end{equation} 
In particular, if $(A,\phi)$ satisfies the vortex equations, then we have the identity 
\[\varepsilon^2 \, |F_A|^2 = \frac{1}{4\varepsilon^2} \, (1 - |\phi|^2).\] 
This relation will play an important role in our subsequent arguments. \\

The equation (1) can be viewed as a simplified version of the Ginzburg-Landau equations (2), 
(3). The Ginzburg-Landau equations play an important role in mathematical physics, where they 
arise in the mathematical description of superconductivity. They have been studied intensively, 
in particular by A. Jaffe and C. H. Taubes \cite{JT,Ta1}. S. Bradlow \cite{Br1,Br2} generalized 
the vortex equations (4), (5) to holomorphic vector bundles over K\"ahler manifolds. \\

The following result shows that every nondegenerate minimal submanifold of $M$ can be obtained 
as the limit of a family of solutions of the Ginzburg-Landau equations. \\

\begin{theorem}
Let $M$ be a Riemannian manifold of dimension $n$, and let $S$ be a nondegenerate minimal 
submanifold of dimension $n - 2$. Then the Ginzburg-Landau equation has a solution for all $0 < 
\varepsilon < \varepsilon_0$. The solutions satisfy 
\[\varepsilon^2 \, |F_A|^2 + |D_A \phi|^2 + \frac{1}{4\varepsilon^2} \, (1 - |\phi|^2)^2 \to 
dH^{n-2}|_S\] 
as $\varepsilon \to 0$. Moreover, the first Chern class of $L$ is the Poincar\'e dual of the 
homology class of $S$.
\end{theorem}

\vspace{2mm}

A related result was established by C. H. Taubes \cite{Ta2,Ta3} for the Seiberg-Witten equations 
on symplectic $4$-manifolds. In this case, every pseudo-holomorphic curve can be approximated by 
a sequence of Seiberg-Witten solutions with parameters $r_j \to \infty$. \\

A similar approach is used in a recent work of F. Pacard and M. Ritor\'e \cite{PR1} which 
relates constant mean curvature hypersurfaces to the theory of phase transitions. \\

Moreover, T. Ilmanen \cite{Il} proved that Brakke's motion by mean curvature is given by the 
limit of a sequence of solutions to the Allen-Cahn equation. \\

In Section 2, we recall some results about the linearized operator on $\mathbb{R}^2$. In 
particular, the kernel of the linearized operator on $\mathbb{R}^2$ is isomorphic to the space 
of parallel vector fields on $\mathbb{R}^2$ (see \cite{JT} and \cite{Ta3} for details). \\

In Section 3, we study the mapping properties of a model operator on the product manifold 
$\mathbb{R}^{n-2} \times \mathbb{R}^2$. \\

In Section 4, we construct a family of approximate solutions of the Ginzburg-Landau equations. 
More precisely, given any normal vector field $v$ satisfying 
\[\|v\|_{\mathcal{C}^{2,\gamma}(S)} \leq \varepsilon,\] 
we construct a pair $(A,\phi)$ such that 
\[\Big \| \Big ( d^* F_A - \frac{1}{2\varepsilon^2} \, (\phi \, \overline{D_A \phi} - \overline
{\phi} \, D_A \phi),D_A^* D_A \phi - \frac{1}{2\varepsilon^2} \, (1 - |\phi|^2) \, \phi \Big ) 
\Big \|_{\mathcal{C}_{\mu,\varepsilon}^\gamma(M)} \leq C\] 
for some $\mu > 0$. Here, the weighted H\"older space $\mathcal{C}_{\mu,\varepsilon}^\gamma(M)$ 
is defined by 
\begin{align*} 
&\|u\|_{\mathcal{C}_{\mu,\varepsilon}^\gamma(M)} \\ 
&= \sup \, e^{\frac{\mu \, \dist(p,S)}{\varepsilon}} \, |u(p)| \\ 
&+ \sup_{\dist(p_1,p_2) \leq \varepsilon} \, \varepsilon^\gamma \, e^{\frac{\mu \, (\dist(p_1,S) 
+ \dist(p_2,S))}{2\varepsilon}} \, \frac{|u(p_1) - u(p_2)|}{\dist(p_1,p_2)^\gamma}. 
\end{align*} 
Moreover, we define 
\[\|(a,f)\|_{\mathcal{C}_{\mu,\varepsilon}^\gamma(M)} = \varepsilon \, \|a\|_{\mathcal{C}_{\mu,
\varepsilon}^\gamma(M)} + \|f\|_{\mathcal{C}_{\mu,\varepsilon}^\gamma(M)}.\]

In Section 5, we derive uniform estimates for the operator $\mathbb{L}_\varepsilon = L
_\varepsilon + T_\varepsilon T_\varepsilon^*$. Here, $L_\varepsilon$ is the linearization of the 
Ginzburg-Landau equations at an approximate solution $(A,\phi)$. Moreover, the operator $T
_\varepsilon$ is defined as 
\[T_\varepsilon u = \big ( \frac{1}{\varepsilon} \, du,-\frac{1}{\varepsilon} \, \phi \, u \big 
)\] 
for $\overline{u} = -u$. Its adjoint is given by 
\[T_\varepsilon^* (a,f) = \varepsilon \, d^* a + \frac{1}{2\varepsilon} \, (\phi \, \overline{f} 
- \overline{\phi} \, f).\] 
The additional term $T_\varepsilon T_\varepsilon^*$ is necessary, because $L_\varepsilon$ is not 
an elliptic operator. \\

To derive uniform estimates independent of $\varepsilon$, we need to restrict the operator 
$\mathbb{L}_\varepsilon$ to a subspace $\mathcal{E}_{\mu,\varepsilon}^\gamma(M) \subset \mathcal
{C}_{\mu,\varepsilon}^\gamma(M)$. A pair $(a,f)$ belongs to $\mathcal{E}_{\mu,\varepsilon}
^\gamma(M)$ if 
\[\int_{NS_x} \varepsilon^2 \, \sum_{\alpha=1}^4 \langle a(e_\alpha^\perp),F_A(w,e_\alpha^\perp) 
\rangle + \int_{NS_x} \langle f,D_{A,w} \phi \rangle = 0\] 
for all $x \in S$ and $w \in NS_x$. \\

In Section 6, we apply the contraction mapping principle to deform the approximate solution 
$(A,\phi)$ to a nearby pair $(\tilde{A},\tilde{\phi})$ such that 
\[(\Id - \mathbb{P}) \Big ( d^* F_{\tilde{A}} - \frac{1}{2\varepsilon^2} \, (\tilde{\phi} \, 
\overline{D_{\tilde{A}} \tilde{\phi}} - \overline{\tilde{\phi}} \, D_{\tilde{A}} \tilde{\phi}) + 
\frac{1}{\varepsilon} \, du,D_{\tilde{A}}^* D_{\tilde{A}} \tilde{\phi} - \frac{1}{2\varepsilon
^2} \, (1 - |\tilde{\phi}|^2) \, \tilde{\phi} - \frac{1}{\varepsilon} \, \tilde{\phi} \, u \Big 
) = 0.\] 
Here, $(\Id - \mathbb{P})$ is the fibrewise projection from $\mathcal{C}_{\mu,\varepsilon}
^\gamma(M)$ to the subspace $\mathcal{E}_{\mu,\varepsilon}^\gamma(M)$. \\

In Section 7, we show that the glueing data can be chosen such that the corresponding pair 
$(\tilde{A},\tilde{\phi})$ satisfies the balancing condition 
\[\mathbb{P} \Big ( d^* F_{\tilde{A}} - \frac{1}{2\varepsilon^2} \, (\tilde{\phi} \, 
\overline{D_{\tilde{A}} \tilde{\phi}} - \overline{\tilde{\phi}} \, D_{\tilde{A}} \tilde{\phi}) + 
\frac{1}{\varepsilon} \, du,D_{\tilde{A}}^* D_{\tilde{A}} \tilde{\phi} - \frac{1}{2\varepsilon
^2} \, (1 - |\tilde{\phi}|^2) \, \tilde{\phi} - \frac{1}{\varepsilon} \, \tilde{\phi} \, u \Big 
) = 0.\] 
This last step of the proof uses the invertibility of the Jacobi operator of $S$. \\

\section{The kernel of the linearized operator on $\mathbb{R}^2$}

In this section, we study the linearized operator on $\mathbb{R}^2$, in particular its kernel. 
To this end, we define an inner product on the space of pairs $(a,f)$ by 
\[\langle (a_1,f_1),(a_2,f_2) \rangle = \int_{\mathbb{R}^2} \varepsilon^2 \, \langle a_1,a_2 
\rangle + \int_{\mathbb{R}^2} \langle f_1,f_2 \rangle.\] 
Then the linearized operator on $\mathbb{R}^2$ satisfies 
\[\langle L_\varepsilon (a,f),(a,f) \rangle = \int_{\mathbb{R}^2} \big | \varepsilon * \! d(ia) 
+ \frac{1}{2\varepsilon} \, (\psi \, \overline{f} + \overline{\psi} \, f) \big |^2 + 2 \int
_{\mathbb{R}^2} \big | \bar{\partial}_B f + \frac{1}{2} \, \psi \, \alpha \big |^2,\] 
where 
\[\alpha = a_1 + ia_2\] 
denotes the $(0,1)$-parts of $a$. This implies 
\begin{align*} 
\partial \alpha &= \frac{1}{2} \, (\partial_1 - i\partial_2) \, (a_1 + ia_2) \\ 
&= \frac{i}{2} \, (\partial_1 a_2 - \partial_2 a_1) + \frac{1}{2} \, (\partial_1 a_1 + \partial
_2 a_2) \\ 
&= \frac{1}{2} \, * \! d(ia) - \frac{1}{2} \, d^* a. 
\end{align*} 
We define an operator $T_\varepsilon: \Omega^0(\mathbb{R}^2,i\mathbb{R}) \to \Omega^1(\mathbb{R}
^2,i\mathbb{R}) \oplus \Omega^0(\mathbb{R}^2,L)$ by 
\[T_\varepsilon u = \big ( \frac{1}{\varepsilon} \, du,-\frac{1}{\varepsilon} \, \psi \, u 
\big )\] 
for $\overline{u} = -u$. Its adjoint is given by 
\[T_\varepsilon^* (a,f) = \varepsilon \, d^* a + \frac{1}{2\varepsilon} \, (\psi \, \overline{f} 
- \overline{\psi} \, f).\] 
This implies 
\[\langle T_\varepsilon T_\varepsilon^* (a,f),(a,f) \rangle = \int_{\mathbb{R}^2} \big | 
\varepsilon \, d^* a + \frac{1}{2\varepsilon} \, (\psi \, \overline{f} - \overline{\psi} \, f) 
\big |^2.\] 
Thus, we conclude that 
\begin{align*} 
&\langle L_\varepsilon (a,f) + T_\varepsilon T_\varepsilon^* (a,f),(a,f) \rangle \\ 
&= \int_{\mathbb{R}^2} \big | \varepsilon * \! d(ia) - \varepsilon \, d^* a + \frac{1}{2
\varepsilon} \, (\psi \, \overline{f} + \overline{\psi} \, f) - \frac{1}{2\varepsilon} \, (\psi 
\, \overline{f} - \overline{\psi} \, f) \big |^2 + 2 \int_{\mathbb{R}^2} \big | \bar{\partial}_B 
f + \frac{1}{2} \, \psi \, \alpha \big |^2 \\ 
&= 4 \int_{\mathbb{R}^2} \big |\varepsilon \, \partial \alpha + \frac{1}{2\varepsilon} \, 
\overline{\psi} \, f \big |^2 + 2 \int_{\mathbb{R}^2} \big | \bar{\partial}_B f + \frac{1}{2} \, 
\psi \, \alpha \big |^2. 
\end{align*}
In particular, the sum 
\[L_\varepsilon + T_\varepsilon T_\varepsilon^*: \Omega^1(\mathbb{R}^2,i\mathbb{R}) \oplus 
\Omega^0(\mathbb{R}^2,L) \to \Omega^1(\mathbb{R}^2,i\mathbb{R}) \oplus \Omega^0(\mathbb{R}^2,L)
\] 
is an elliptic operator. \\

\begin{proposition}
The kernel of the operator $L_\varepsilon + T_\varepsilon T_\varepsilon^*$ is a vector space of 
real dimension $2$. It consists of all pairs of the form $(F_B(w,\cdot),D_{B,w} \psi)$, where 
$w$ is a fixed vector in $\mathbb{R}^2$.
\end{proposition}

\textbf{Proof.} 
Suppose that $(a,f)$ satisfies the equation 
\[L_\varepsilon (a,f) + T_\varepsilon T_\varepsilon^* (a,f) = 0.\] 
This implies 
\[\varepsilon \, \partial \alpha + \frac{1}{2\varepsilon} \, \overline{\psi} \, f = 0\] 
and 
\[\bar{\partial}_B f + \frac{1}{2} \, \psi \, \alpha = 0.\] 
The set $V$ of pairs $(a,f)$ satisfying these conditions is a vector space of real dimension $2$ 
(see, for example, \cite{JT}). We claim that the pair 
\[a = F_B(w,\cdot)\] 
and 
\[f = D_{B,w} \psi\] 
belongs to $V$. Using the identity 
\[\varepsilon * \! (iF_B) = \frac{1}{2\varepsilon} \, (1 - |\psi|^2),\] 
we obtain 
\begin{align*} 
\varepsilon * \! d(ia) + \frac{1}{2\varepsilon} \, (\psi \, \overline{f} + \overline{\psi} \, f) 
&= \varepsilon * \! d(iF_B(w,\cdot)) + \frac{1}{2\varepsilon} \, (\psi \, \overline{D_{B,w} 
\psi} + \overline{\psi} \, D_{B,w} \psi) \\ 
&= \varepsilon \, \partial_w * \! (iF_B) + \frac{1}{2\varepsilon} \, \partial_w |\psi|^2 \\ 
&= 0 
\end{align*}
and 
\begin{align*} 
\varepsilon \, d^* a + \frac{1}{2\varepsilon} \, (\psi \, \overline{f} - \overline{\psi} \, f) 
&= -\varepsilon \, (d^* F_B)(w) + \frac{1}{2\varepsilon} \, (\psi \, \overline{D_{B,w} \psi} - 
\overline{\psi} \, D_{B,w} \psi) \\ 
&= 0. 
\end{align*} 
Moreover, we have 
\[\bar{\partial}_B f = \bar{\partial}_B D_{B,w} \psi = \bar{\partial}_B D_{B,w} \psi - D_{B,w} 
\bar{\partial}_B \psi = -\frac{1}{2} \, \alpha \, \psi.\] 
This proves the assertion. \\

\section{The model problem on $\mathbb{R}^{n-2} \times \mathbb{R}^2$}

In this section, we consider a complex line bundle $L$ over the product $\mathbb{R}^{n-2} \times 
\mathbb{R}^2$. Let $B$ be a connection on $L$ and let $\psi$ be a section of $L$. We assume that 
the pair $(B,\psi)$ is invariant under translations along the $\mathbb{R}^{n-2}$ factor and 
agrees with the one-vortex solution along the $\mathbb{R}^2$ factor. \\

As in Section 2, we consider the inner product 
\[\langle (a_1,f_1),(a_2,f_2) \rangle = \int_{\mathbb{R}^{n-2} \times \mathbb{R}^2} \varepsilon
^2 \, \langle a_1,a_2 \rangle + \int_{\mathbb{R}^{n-2} \times \mathbb{R}^2} \langle f_1,f_2 
\rangle.\] 
Then the linearized operator $L_\varepsilon$ satisfies 
\begin{align*} 
\langle L_\varepsilon(a,f),(a,f) \rangle 
&= \int_{\mathbb{R}^{n-2} \times \mathbb{R}^2} (\varepsilon^2 \, |da|^2 + |D_B f|^2) \\ 
&+ 2 \int_{\mathbb{R}^{n-2} \times \mathbb{R}^2} (\langle D_B \psi,a \, f \rangle + \langle a \, 
\psi,D_B f \rangle) \\ 
&+ \int_{\mathbb{R}^{n-2} \times \mathbb{R}^2} \Big ( |\psi|^2 \, |a|^2 + \frac{1}{\varepsilon
^2} \, \Re(\overline{\psi} \, f)^2 - \frac{1}{2\varepsilon^2} \, (1 - |\psi|^2) \, |f|^2 \Big ). 
\end{align*} 
As in Section 2, we define an operator $T_\varepsilon: \Omega^0(\mathbb{R}^{n-2} \times \mathbb
{R}^2,i\mathbb{R}) \to \Omega^1(\mathbb{R}^{n-2} \times \mathbb{R}^2,i\mathbb{R}) \oplus \Omega
^0(\mathbb{R}^{n-2} \times \mathbb{R}^2,L)$ by 
\[T_\varepsilon u = \big ( \frac{1}{\varepsilon} \, du,-\frac{1}{\varepsilon} \, \psi \, u \big 
)\] 
for $\overline{u} = -u$. Its adjoint is given by 
\[T_\varepsilon^* (a,f) = \varepsilon \, d^* a + \frac{1}{2\varepsilon} \, (\psi \, \overline{f} 
- \overline{\psi} \, f).\] 
This implies 
\[\langle T_\varepsilon T_\varepsilon^* (a,f),(a,f) \rangle = \int_{\mathbb{R}^{n-2} \times 
\mathbb{R}^2} \Big | \varepsilon \, d^* a + \frac{1}{2\varepsilon} \, (\psi \, \overline{f} - 
\overline{\psi} f) \Big |^2.\] 
Therefore, we obtain 
\begin{align*} 
&\langle L_\varepsilon (a,f),(a,f) \rangle + \langle T_\varepsilon T_\varepsilon^* (a,f),(a,f) 
\rangle \\ 
&= \int_{\mathbb{R}^{n-2} \times \mathbb{R}^2} (\varepsilon^2 \, |da|^2 + \varepsilon^2 \, |d^* 
a|^2 + |D_B f|^2) \\ 
&+ 2 \int_{\mathbb{R}^{n-2} \times \mathbb{R}^2} (\langle D_B \psi,a \, f \rangle + \langle a \, 
\psi,D_B f \rangle) - 2 \int_{\mathbb{R}^{n-2} \times \mathbb{R}^2} \langle d^* a \, \psi,f 
\rangle \\ 
&+ \int_{\mathbb{R}^{n-2} \times \mathbb{R}^2} \Big ( |\psi|^2 \, |a|^2 + \frac{1}{\varepsilon
^2} \, \Re(\overline{\psi} \, f)^2 + \frac{1}{\varepsilon^2} \, \Im(\overline{\psi} \, f)^2 - 
\frac{1}{2\varepsilon^2} \, (1 - |\psi|^2) \, |f|^2 \Big ) \\ 
&= \int_{\mathbb{R}^{n-2} \times \mathbb{R}^2} (\varepsilon^2 \, |da|^2 + \varepsilon^2 \, |d^* 
a|^2 + |D_B f|^2) \\ 
&+ 4 \int_{\mathbb{R}^{n-2} \times \mathbb{R}^2} \langle D_B \psi,a \, f \rangle \\ 
&+ \int_{\mathbb{R}^{n-2} \times \mathbb{R}^2} \Big ( |\psi|^2 \, |a|^2 + \frac{1}{\varepsilon
^2} \, |\psi|^2 \, |f|^2 - \frac{1}{2\varepsilon^2} \, (1 - |\psi|^2) \, |f|^2 \Big ). 
\end{align*} 
From this it follows that 
\begin{align*} 
L_\varepsilon (a,f) + T_\varepsilon T_\varepsilon^* (a,f) = \Big ( 
&\nabla^* \nabla a - \frac{1}{\varepsilon^2} \, (\overline{D_B \psi} \, f - D_B \psi \, 
\overline{f}) + \frac{1}{\varepsilon^2} \, |\psi|^2 \, a, \\ &D_B^* D_B f - 2 \, *(a \wedge * 
D_B \psi) + \frac{1}{\varepsilon^2} \, |\psi|^2 \, f - \frac{1}{2\varepsilon^2} \, (1 - |\psi|
^2) \, f \Big ) 
\end{align*} 
for $a \in \Omega^1(\mathbb{R}^{n-2} \times \mathbb{R}^2,i\mathbb{R})$ and 
$f \in \Omega^0(\mathbb{R}^{n-2} \times \mathbb{R}^2,L)$. For abbreviation, let $\mathbb{L}
_\varepsilon = L_\varepsilon + T_\varepsilon T_\varepsilon^*$. Note that 
\begin{align*} 
\mathbb{L}_\varepsilon: \; &\Omega^1(\mathbb{R}^{n-2} \times \mathbb{R}^2,i\mathbb{R}) \oplus 
\Omega^0(\mathbb{R}^{n-2} \times \mathbb{R}^2,L) \\ 
&\to \Omega^1(\mathbb{R}^{n-2} \times \mathbb{R}^2,i\mathbb{R}) \oplus \Omega^0(\mathbb{R}^{n-2} 
\times \mathbb{R}^2,L) 
\end{align*}
is an elliptic operator. \\

We define the weighted H\"older space $\mathcal{C}_{\mu,\varepsilon}^\gamma(\mathbb{R}^{n-2} 
\times \mathbb{R}^2)$ by 
\begin{align*} 
&\|u\|_{\mathcal{C}_{\mu,\varepsilon}^\gamma(\mathbb{R}^{n-2} \times \mathbb{R}^2)} \\ 
&= \sup \, e^{\frac{\mu \, |y|}{\varepsilon}} \, |u(x,y)| \\ 
&+ \sup_{|x_1 - x_2| + |y_1 - y_2| \leq \varepsilon} \, \varepsilon^\gamma \, e^{\frac{\mu \, 
(|y_1| + |y_2|)}{2\varepsilon}} \, \frac{|u(x_1,y_1) - u(x_2,y_2)|}{(|x_1 - x_2| + |y_1 - y_2|)
^\gamma}. 
\end{align*} 
More generally, let 
\[\|u\|_{\mathcal{C}_{\mu,\varepsilon}^{k,\gamma}(\mathbb{R}^{n-2} \times \mathbb{R}^2)} = \sum
_{l=0}^k \varepsilon^l \, \|\nabla^l u\|_{\mathcal{C}_{\mu,\varepsilon}^\gamma(\mathbb{R}^{n-2} 
\times \mathbb{R}^2)}.\] 
Furthermore, we define 
\[\|(a,f)\|_{\mathcal{C}_{\mu,\varepsilon}^{k,\gamma}(\mathbb{R}^{n-2} \times \mathbb{R}^2)} = 
\varepsilon \, \|a\|_{\mathcal{C}_{\mu,\varepsilon}^{k,\gamma}(\mathbb{R}^{n-2} \times \mathbb
{R}^2)} + \|f\|_{\mathcal{C}_{\mu,\varepsilon}^{k,\gamma}(\mathbb{R}^{n-2} \times \mathbb{R}
^2)}.\]

Let $\mathcal{E}_{\mu,\varepsilon}^{k,\gamma}(\mathbb{R}^{n-2} \times \mathbb{R}^2)$ be the set 
of all pairs $(a,f) \in \Omega^1(\mathbb{R}^{n-2} \times \mathbb{R}^2,i\mathbb{R}) \oplus \Omega
^0(\mathbb{R}^{n-2} \times \mathbb{R}^2,L)$ such that $(a,f) \in \mathcal{C}_{\mu,\varepsilon}
^{k,\gamma}(\mathbb{R}^{n-2} \times \mathbb{R}^2)$ and 
\[\int_{\{x\} \times \mathbb{R}^2} \varepsilon^2 \, \sum_{\alpha=1}^2 \langle a(e_\alpha^\perp),
F_B(w,e_\alpha^\perp) \rangle + \int_{\{x\} \times \mathbb{R}^2} \langle f,D_{B,w} \psi \rangle 
= 0\] 
for all $x \in \mathbb{R}^{n-2}$ and all $w \in \mathbb{R}^2$. \\

\begin{proposition}
The operator $\mathbb{L}_\varepsilon$ maps $\mathcal{E}_{\mu,\varepsilon}^{2,\gamma}(\mathbb{R}
^{n-2} \times \mathbb{R}^2)$ into $\mathcal{E}_{\mu,\varepsilon}^\gamma(\mathbb{R}^{n-2} \times 
\mathbb{R}^2)$.
\end{proposition}

\textbf{Proof.} 
It is obvious from the definition that $\mathbb{L}_\varepsilon$ maps $\mathcal{C}_{\mu,
\varepsilon}^{2,\gamma}(\mathbb{R}^{n-2} \times \mathbb{R}^2)$ into $\mathcal{C}_{\mu,
\varepsilon}^\gamma(\mathbb{R}^{n-2} \times \mathbb{R}^2)$. We now assume that $(a,f) \in 
\mathcal{C}_{\mu,\varepsilon}^{2,\gamma}(\mathbb{R}^{n-2} \times \mathbb{R}^2)$ satisfies 
\[\int_{\mathbb{R}^2} \varepsilon^2 \langle a,F_B(w,\cdot) \rangle + \int_{\mathbb{R}^2} \langle 
f,D_{B,w} \psi \rangle = 0\] 
for all $x \in \mathbb{R}^{n-2}$ and all $w \in \mathbb{R}^2$. Taking derivatives in horizontal 
direction, we obtain 
\[\int_{\{x\} \times \mathbb{R}^2} \varepsilon^2 \, \sum_{\alpha=1}^2 \sum_{j=1}^{n-2} \langle 
\partial_j \partial_j a(e_\alpha^\perp),F_B(w,e_\alpha^\perp) \rangle + \int_{\mathbb{R}^2} 
\langle f,D_{B,w} \psi \rangle = 0.\] 
Furthermore, integration by parts gives 
\begin{align*} 
&\int_{\{x\} \times \mathbb{R}^2} \sum_{\alpha=1}^2 \sum_{\rho=1}^2 \varepsilon^2 \, \langle 
\nabla_{e_\rho^\perp} \nabla_{e_\rho^\perp} a(e_\alpha^\perp),F_B(w,e_\alpha^\perp) \rangle \\ 
&+ \int_{\{x\} \times \mathbb{R}^2} \sum_{\rho=1}^2 \langle D_{B,e_\rho^\perp} D_{B,e_\rho
^\perp} f,D_{B,w} \psi \rangle \\ 
&- 2 \int_{\{x\} \times \mathbb{R}^2} \langle D_B \psi,a \, D_{B,w} \psi \rangle - 2 \int_{\{x\} 
\times \mathbb{R}^2} \langle D_B \psi,F_B(w,\cdot) \, f \rangle \\ 
&- \int_{\{x\} \times \mathbb{R}^2} \Big ( |\psi|^2 \, \langle a,F_B(w,\cdot) \rangle + \frac
{1}{\varepsilon^2} \, |\psi|^2 \, \langle f,D_{B,w} \psi \rangle - \frac{1}{2\varepsilon^2} \, 
(1 - |\psi|^2) \, \langle f,D_{B,w} \psi \rangle \Big ) \\ 
&= \int_{\{x\} \times \mathbb{R}^2} \sum_{\alpha=1}^2 \sum_{\rho=1}^2 \varepsilon^2 \, \langle 
a(e_\alpha^\perp),\nabla_{e_\rho^\perp} \nabla_{e_\rho^\perp} F_B(w,e_\alpha^\perp) \rangle \\ 
&+ \int_{\{x\} \times \mathbb{R}^2} \sum_{\rho=1}^2 \langle f,D_{B,e_\rho^\perp} D_{B,e_\rho
^\perp} D_{B,w} \psi \rangle \\ 
&- 2 \int_{\{x\} \times \mathbb{R}^2} \langle D_B \psi,a \, D_{B,w} \psi \rangle - 2 \int_{\{x\} 
\times \mathbb{R}^2} \langle D_B \psi,F_B(w,\cdot) \, f \rangle \\ 
&- \int_{\{x\} \times \mathbb{R}^2} \Big ( |\psi|^2 \, \langle a,F_B(w,\cdot) \rangle + \frac
{1}{\varepsilon^2} \, |\psi|^2 \, \langle f,D_{B,w} \psi \rangle - \frac{1}{2\varepsilon^2} \, 
(1 - |\psi|^2) \, \langle f,D_{B,w} \psi \rangle \Big ) \\ 
&= 0.
\end{align*} 
Hence, if we define $(b,h) = \mathbb{L}_\varepsilon (a,f)$, then we obtain 
\[\int_{\{x\} \times \mathbb{R}^2} \varepsilon^2 \, \sum_{\alpha=1}^2 \langle b(e_\alpha^\perp),
F_B(w,e_\alpha^\perp) \rangle + \int_{\{x\} \times \mathbb{R}^2} \langle h,D_{B,w} \psi \rangle 
= 0\] 
for all $x \in \mathbb{R}^{n-2}$ and $w \in \mathbb{R}^2$. \\

\begin{proposition}
Let $0 < \nu < 1$, $(b,h) \in \mathcal{C}_{\mu,\varepsilon}^\gamma(\mathbb{R}^2)$, and $\eta \in 
\mathcal{S}(\mathbb{R}^{n-2})$. Moreover, assume that the Fourier transform of $\eta$ satisfies 
$\hat{\eta}(\xi) = 0$ for $|\xi| \leq \delta$ for some $\delta > 0$. Then there exists a pair 
$(a,f) \in \mathcal{C}_{\mu,\varepsilon}^{2,\gamma}(\mathbb{R}^{n-2} \times \mathbb{R}^2)$ such 
that 
\[\mathbb{L}_\varepsilon (a,f) = (\eta(x) \, b(y),\eta(x) \, h(y)).\]
\end{proposition}

\textbf{Proof.} 
We perform a Fourier transformation in the $\mathbb{R}^{n-2}$ variables. Let 
\[\eta(x) = \int_{\mathbb{R}^{n-2}} e^{ix\xi} \, \hat{\eta}(\xi) \, d\xi.\] 
For every $\xi \in \mathbb{R}^{n-2}$, there exists a pair $(\hat{a}(\xi,\cdot),\hat{f}(\xi,
\cdot)) \in \mathcal{C}_{\mu,\varepsilon}^{2,\gamma}(\mathbb{R}^2)$ such that 
\begin{align*} 
&\sum_{\rho=1}^2 \partial_\rho \partial_\rho \hat{a}(\xi,y) + \frac{1}{\varepsilon^2} \, 
(\overline{D_B \psi(y)} \, \hat{f}(\xi,y) - D_B \psi(y) \, \overline{\hat{f}(\xi,y)}) \\ 
&- \frac{1}{\varepsilon^2} \, |\psi(y)|^2 \, \hat{a}(\xi,y) - |\xi|^2 \, \hat{a}(\xi,y) = -b(y) 
\end{align*} 
and 
\begin{align*} 
&\sum_{\rho=1}^2 D_{B,\rho} D_{B,\rho} \hat{f}(\xi,y) + 2 \sum_{\rho=1}^2 D_{B,\rho} \psi(y) \, 
\hat{a}_\rho(\xi,y) \\ &- \frac{1}{\varepsilon^2} \, |\psi(y)|^2 \, \hat{f}(\xi,y) + \frac{1}{2
\varepsilon^2} \, (1 - |\psi(y)|^2) \, \hat{f}(\xi,y) - |\xi|^2 \, \hat{f}(\xi,y) = -h(y). 
\end{align*}
We now define a pair $(a,f)$ by 
\[a(x,y) = \int_{\mathbb{R}^{n-2}} e^{ix\xi} \, \hat{\eta}(\xi) \, \hat{a}(\xi,y) \, d\xi\] 
and 
\[f(x,y) = \int_{\mathbb{R}^{n-2}} e^{ix\xi} \, \hat{\eta}(\xi) \, \hat{f}(\xi,y) \, d\xi.\] 
Then the pair $(a,f)$ satisfies 
\[\sum_{i=1}^{n-2} \partial_i \partial_i a + \sum_{\rho=1}^2 \partial_\rho \partial_\rho a + 
\frac{1}{\varepsilon^2} \, (\overline{D_B \psi} \, f - D_B \psi \, \overline{f}) - \frac
{1}{\varepsilon^2} \, |\psi|^2 \, a = -\eta(x) \, b(y)\] 
and 
\[\sum_{i=1}^{n-2} \partial_i \partial_i f + \sum_{\rho=1}^2 D_{B,\rho} D_{B,\rho} f + 2 \sum
_{\rho=1}^2 D_{B,\rho} \psi \, a_\rho - \frac{1}{\varepsilon^2} \, |\psi|^2 \, f + \frac{1}{2
\varepsilon^2} \, (1 - |\psi|^2) \, \hat{f} = -\eta(x) \, h(y).\] 
From this we deduce that 
\[\mathbb{L}_\varepsilon (a,f) = (\eta(x) \, b(y),\eta(x) \, h(y)).\] 
This proves the assertion. \\

\begin{proposition}
Let $0 < \nu < 1$, and suppose that $(a,f) \in \mathcal{E}_{\mu,\varepsilon}^{2,\gamma}(\mathbb
{R}^{n-2} \times \mathbb{R}^2)$ satisfies $\mathbb{L}_\varepsilon (a,f) = 0$. Then $(a,f) = 0$.
\end{proposition}

\textbf{Proof.} 
Let $(b,h) \in \mathcal{C}_{\mu,\varepsilon}^\gamma(\mathbb{R}^2)$ and $\zeta \in \mathcal{S}
(\mathbb{R}^{n-2})$ be given. We define a function $\eta \in \mathcal{S}(\mathbb{R}^{n-2})$ by 
$\eta(x) = \zeta(x + x_0) - \zeta(x)$. Then the Fourier transform of $\eta$ satisfies $\hat
{\eta}(0) = 0$. We approximate $\eta$ by functions $\eta_\delta$ such that 
\[\hat{\eta}_\delta(\xi) = 0\] 
for $|\xi| \leq \delta$ and 
\[\hat{\eta}_\delta(\xi) = \hat{\eta}(\xi)\] 
for $|\xi| \geq 2\delta$. Using the condition $\hat{\eta}(0) = 0$, we obtain 
\[\|D^{n-2}(\hat{\eta} - \hat{\eta}_\delta)\|_{L^{\frac{p}{p-1}}(\mathbb{R}^{n-2})} \leq C \, 
\delta^{1-\frac{n-2}{p}},\] 
hence 
\[\big \| (1 + |x|)^{n-2} \, (\eta - \eta_\delta) \big \|_{L^p(\mathbb{R}^{n-2})} \leq C \, 
\delta^{1-\frac{n-2}{p}}\] 
for all $p \geq 2$. From this it follows that 
\begin{align*} 
\|\eta - \eta_\delta\|_{L^1(\mathbb{R}^{n-2})} 
&\leq \big \| (1 + |x|)^{-(n-2)} \big \|_{L^{\frac{p}{p-1}}(\mathbb{R}^{n-2})} \, \big \| (1 + 
|x|)^{n-2} \, (\eta - \eta_\delta) \big \|_{L^p(\mathbb{R}^{n-2})} \\ 
&\leq C \, \delta^{1-\frac{n-2}{p}} 
\end{align*} 
for all $p \geq 2$. This implies 
\[\|\eta - \eta_\delta\|_{L^1(\mathbb{R}^{n-2})} \to 0\] 
as $\delta \to 0$. \\

For each $\delta > 0$, the pair $(\eta_\delta(x) \, b(y),\eta_\delta(x) \, h(y))$ belongs to the 
image of $\mathbb{L}_\varepsilon$. Since $(a,f)$ belongs to the kernel of $\mathbb{L}
_\varepsilon$, we obtain 
\[\int_{\mathbb{R}^{n-2} \times \mathbb{R}^2} \varepsilon^2 \, \langle a(x,y),\eta_\delta(x) \, 
b(y) \rangle + \int_{\mathbb{R}^{n-2} \times \mathbb{R}^2} \langle f(x,y),\eta_\delta(x) \, h(y) 
\rangle = 0.\] 
Letting $\delta \to 0$, we obtain 
\[\int_{\mathbb{R}^{n-2} \times \mathbb{R}^2} \varepsilon^2 \, \langle a(x,y),\eta(x) \, b(y) 
\rangle + \int_{\mathbb{R}^{n-2} \times \mathbb{R}^2} \langle f(x,y),\eta(x) \, h(y) \rangle = 
0.\] 
From this it follows that 
\begin{align*} 
&\int_{\mathbb{R}^{n-2} \times \mathbb{R}^2} \varepsilon^2 \, \langle a(x,y),\zeta(x) \, b(y) 
\rangle + \int_{\mathbb{R}^{n-2} \times \mathbb{R}^2} \langle f(x,y),\zeta(x) \, h(y) \rangle \\ 
&= \int_{\mathbb{R}^{n-2} \times \mathbb{R}^2} \varepsilon^2 \, \langle a(x-x_0,y),\zeta(x) \, 
b(y) \rangle + \int_{\mathbb{R}^{n-2} \times \mathbb{R}^2} \langle f(x-x_0,y),\zeta(x) \, h(y) 
\rangle. 
\end{align*} 
Since $(b,h)$ and $\zeta$ are arbitrary, we conclude that $a(x,y) = a(x-x_0,y)$ and $f(x,y) = 
f(x-x_0,y)$. Therefore, $a(x,y)$ and $f(x,y)$ are constant in $x$. Using Proposition 2.1, we 
obtain 
\[(a,f) = (F_B(w,\cdot),D_{B,w} \psi)\] 
for some $w \in \mathbb{R}^2$. This proves the assertion. \\

\begin{proposition}
Let $0 < \nu < 1$. Then we have the estimate 
\[\|(a,f)\|_{\mathcal{C}_{\mu,\varepsilon}^{2,\gamma}(\mathbb{R}^{n-2} \times \mathbb{R}^2)} 
\leq C \, \varepsilon^2 \, \|\mathbb{L}_\varepsilon (a,f)\|_{\mathcal{C}_{\mu,\varepsilon}
^\gamma(\mathbb{R}^{n-2} \times \mathbb{R}^2)}\] 
for all $(a,f) \in \mathcal{E}_{\mu,\varepsilon}^{2,\gamma}(\mathbb{R}^{n-2} \times \mathbb{R}
^2)$.
\end{proposition}

\textbf{Proof.} 
By Schauder estimates, it suffices to show that 
\[\sup \, e^{\frac{\mu |y|}{\varepsilon}} \, (\varepsilon \, |a(x,y)| + |f(x,y)|) \leq C \, \sup 
\, \varepsilon^2 \, e^{\frac{\mu |y|}{\varepsilon}} \, (\varepsilon \, |b(x,y)| + |h(x,y)|),\] 
where $(b,h) = \mathbb{L}_\varepsilon (a,f)$. To prove this estimate, we argue by contradiction. 
Let $(a^{(j)},f^{(j)})$ be a sequence of pairs such that 
\[\sup \, e^{\frac{\mu |y|}{\varepsilon}} \, (\varepsilon \, |a^{(j)}(x,y)| + |f^{(j)}(x,y)|) = 
1\] 
and 
\[\sup \, \varepsilon^2 \, e^{\frac{\mu |y|}{\varepsilon}} \, (\varepsilon \, |b^{(j)}(x,y)| + 
|h^{(j)}(x,y)|) \to 0,\] 
where $(b^{(j)},h^{(j)}) = \mathbb{L}_{\varepsilon} (a^{(j)},f^{(j)})$. We choose a sequence of 
points $(x_j,y_j) \in \mathbb{R}^{n-2} \times \mathbb{R}^2$ such that 
\[e^{\frac{\mu |y_j|}{\varepsilon}} \, (\varepsilon \, |a^{(j)}(x_j,y_j)| + |f^{(j)}(x_j,y_j)|) 
\geq \frac{1}{2}\] 
for all $j$. There are two possibilities: \\

(i) Suppose that the sequence $|y_j|$ is bounded. After passing to a subsequence, we may assume 
that the sequence $(a^{(j)},f^{(j)})$ converges to a pair $(a,f) \in \mathcal{E}_{\mu,
\varepsilon}^{2,\gamma}(\mathbb{R}^{n-2} \times \mathbb{R}^2)$ such that 
\[\sup \, e^{\frac{\mu |y|}{\varepsilon}} \, (\varepsilon \, |a(x,y)| + |f(x,y)|) \leq 1\] 
and 
\[\mathbb{L}_\varepsilon (a,f) = 0.\] 
Using Proposition 3.3, we conclude that $(a,f) = 0$. This is a contradiction. \\

(ii) We now assume that $|y_j| \to \infty$. We define a sequence of pairs $(\tilde{a}_j,\tilde
{f}_j)$ by \[\tilde{a}^{(j)}(x,y) = e^{\frac{\mu |y_j|}{\varepsilon}} \, a^{(j)}(x + x_j,y + y
_j)\] 
and 
\[\tilde{f}^{(j)}(x,y) = e^{\frac{\mu |y_j|}{\varepsilon}} \, f^{(j)}(x + x_j,y + y_j).\] 
After passing to a subsequence, we may assume that the sequence $(\tilde{a}_j,\tilde{f}_j)$ 
converges to a pair $(\tilde{a},\tilde{f})$. The pair $(\tilde{a},\tilde{f})$ is defined on 
$\mathbb{R}^{n-2} \times \mathbb{R}^2$ and satisfies 
\[\sup \, e^{-\frac{\mu |y|}{\varepsilon}} \, (\varepsilon \, |\tilde{a}(x,y)| + |\tilde{f}(x,
y)|) \leq 1\] 
and 
\[\Big ( \nabla^* \nabla \tilde{a} + \frac{1}{\varepsilon^2} \, \tilde{a},\nabla^* \nabla \tilde
{f} + \frac{1}{\varepsilon^2} \, \tilde{f} \Big ) = 0.\] 
If $\mu$ is sufficiently small, it follows that $(\tilde{a},\tilde{f}) = 0$. This is a 
contradiction. \\

\begin{proposition}
Let $0 < \nu < 1$. Assume that $(b,h) \in \mathcal{C}_{\mu,\varepsilon}^\gamma(\mathbb{R}^2)$ 
satisfies 
\[\int_{\mathbb{R}^2} \varepsilon^2 \, \langle b,F_B(w,\cdot) \rangle + \int_{\mathbb{R}^2} 
\langle h,D_{B,w} \psi \rangle = 0\] 
for all $x \in \mathbb{R}^{n-2}$ and all $w \in \mathbb{R}^2$. Moreover, let $\eta \in \mathcal
{S}(\mathbb{R}^{n-2})$. Then there exists a pair $(a,f) \in \mathcal{E}_{\mu,\varepsilon}^{2,
\gamma}(\mathbb{R}^{n-2} \times \mathbb{R}^2)$ such that 
\[\mathbb{L}_\varepsilon (a,f) = (\eta(x) \, b(y),\eta(x) \, h(y)).\]
\end{proposition}

\textbf{Proof.} 
Let 
\[\eta(x) = \int_{\mathbb{R}^{n-2}} e^{ix\xi} \, \hat{\eta}(\xi) \, d\xi.\] 
For every $\xi \in \mathbb{R}^{n-2}$, there exists a $1$-form $\hat{a}(\xi,\cdot) \in \mathcal
{C}_{\mu,\varepsilon}^{2,\gamma}(\mathbb{R}^2)$ such that 
\begin{align*} 
&\sum_{\rho=1}^2 \partial_\rho \partial_\rho \hat{a}(\xi,y) + \frac{1}{\varepsilon^2} \, 
(\overline{D_B \psi(y)} \, \hat{f}(\xi,y) - D_B \psi(y) \, \overline{\hat{f}(\xi,y)}) \\ 
&- \frac{1}{\varepsilon^2} \, |\psi(y)|^2 \, \hat{a}(\xi,y) - |\xi|^2 \, \hat{a}(\xi,y) = -b(y) 
\end{align*}
and 
\begin{align*} 
&\sum_{\rho=1}^2 D_{B,\rho} D_{B,\rho} \hat{f}(\xi,y) + 2 \sum_{\rho=1}^2 D_{B,\rho} \psi(y) \, 
\hat{a}_\rho(\xi,y) \\ &- \frac{1}{\varepsilon^2} \, |\psi(y)|^2 \, \hat{f}(\xi,y) + \frac{1}{2
\varepsilon^2} \, (1 - |\psi(y)|^2) \, \hat{f}(\xi,y) - |\xi|^2 \, \hat{f}(\xi,y) = -h(y). 
\end{align*}
Furthermore, the pair $(\hat{a}(\xi,\cdot),\hat{f}(\xi,\cdot))$ satisfies 
\[\int_{\mathbb{R}^2} \varepsilon^2 \, \langle \hat{a}(\xi,\cdot),F_B(w,\cdot) \rangle + \int
_{\mathbb{R}^2} \langle \hat{f}(\xi,\cdot),D_{B,w} \psi \rangle = 0\] 
for all $x \in \mathbb{R}^{n-2}$ and all $w \in \mathbb{R}^2$. We now define a pair $(a,f) \in 
\mathcal{E}_{\mu,\varepsilon}^{2,\gamma}(\mathbb{R}^{n-2} \times \mathbb{R}^2)$ by 
\[a(x,y) = \int_{\mathbb{R}^{n-2}} e^{ix\xi} \, \hat{\eta}(\xi) \, \hat{a}(\xi,y) \, d\xi\] 
and 
\[f(x,y) = \int_{\mathbb{R}^{n-2}} e^{ix\xi} \, \hat{\eta}(\xi) \, \hat{f}(\xi,y) \, d\xi.\] 
Then the pair $(a,f)$ satisfies 
\[\sum_{i=1}^{n-2} \partial_i \partial_i a + \sum_{\rho=1}^2 \partial_\rho \partial_\rho a + 
\frac{1}{\varepsilon^2} \, (\overline{D_B \psi} \, f - D_B \psi \, \overline{f}) - \frac
{1}{\varepsilon^2} \, |\psi|^2 \, a = -\eta(x) \, b(y)\] 
and 
\[\sum_{i=1}^{n-2} \partial_i \partial_i f + \sum_{\rho=1}^2 D_{B,\rho} D_{B,\rho} f + 2 \sum
_{\rho=1}^2 D_{B,\rho} \psi \, a_\rho - \frac{1}{\varepsilon^2} \, |\psi|^2 \, f + \frac{1}{2
\varepsilon^2} \, (1 - |\psi|^2) \, \hat{f} = -\eta(x) \, h(y).\] 
Thus, we conclude that 
\[\mathbb{L}_\varepsilon (a,f) = (\eta(x) \, b(y),\eta(x) \, h(y)).\] 
This proves the assertion. \\

\begin{corollary}
Let $0 < \nu < 1$, and suppose that $(b,h) \in \mathcal{E}_{\mu,\varepsilon}^\gamma(\mathbb{R}
^{n-2} \times \mathbb{R}^2)$ has compact support. Then there exists a pair $(a,f) \in \mathcal
{E}_{\mu,\varepsilon}^{2,\gamma}(\mathbb{R}^{n-2} \times \mathbb{R}^2)$ such that 
\[\|(a,f)\|_{\mathcal{C}_{\mu,\varepsilon}^{2,\gamma}(\mathbb{R}^{n-2} \times \mathbb{R}^2)} 
\leq C \, \varepsilon^2 \, \|(b,h)\|_{\mathcal{C}_{\mu,\varepsilon}^\gamma(\mathbb{R}^{n-2} 
\times \mathbb{R}^2)}\] 
and 
\[\mathbb{L}_\varepsilon (a,f) = (b,h).\] 
\end{corollary}

\textbf{Proof.} 
It follows from Proposition 3.4 that the range of the operator $\mathbb{L}_\varepsilon: \mathcal
{E}_{\mu,\varepsilon}^{2,\gamma}(\mathbb{R}^{n-2} \times \mathbb{R}^2) \to \mathcal{E}_{\mu,
\varepsilon}^\gamma(\mathbb{R}^{n-2} \times \mathbb{R}^2)$ is a closed subspace of the Banach 
space $\mathcal{E}_{\mu,\varepsilon}^\gamma(\mathbb{R}^{n-2} \times \mathbb{R}^2)$. By 
Proposition 3.5, it contains all pairs of the form $(\eta(x) \, b(y),\eta(x) \, h(y))$, where 
$\eta \in \mathcal{S}(\mathbb{R}^{n-2})$ and $(b,h) \in \mathcal{C}_{\mu,\varepsilon}^\gamma
(\mathbb{R}^2)$ satisfies 
\[\int_{\mathbb{R}^2} \varepsilon^2 \, \langle b,F_B(w,\cdot) \rangle + \int_{\mathbb{R}^2} 
\langle h,D_{B,w} \psi \rangle = 0\] 
for all $x \in \mathbb{R}^{n-2}$ and all $w \in \mathbb{R}^2$. The assertion follows now by 
approximation. \\

\begin{proposition}
Let $0 < \nu < 1$. Suppose that $(b,h) \in \mathcal{E}_{\mu,\varepsilon}^\gamma(\mathbb{R}^{n-2} 
\times \mathbb{R}^2)$ is supported in the set $\{(x,y) \in \mathbb{R}^{n-2} \times \mathbb{R}^2: 
|x| \leq \delta, \, |y| \leq 2\delta\}$. Then there exists a pair $(a,f) \in \mathcal{C}_{\mu,
\varepsilon}^{2,\gamma}(\mathbb{R}^{n-2} \times \mathbb{R}^2)$ such that $(a,f)$ is supported in 
$\{(x,y) \in \mathbb{R}^{n-2} \times \mathbb{R}^2: |x| \leq 2\delta, \, |y| \leq 4\delta\}$, 
\[\|(a,f)\|_{\mathcal{C}_{\mu,\varepsilon}^{2,\gamma}(\mathbb{R}^{n-2} \times \mathbb{R}^2)} 
\leq C \, \varepsilon^2 \, \|(b,h)\|_{\mathcal{C}_{\mu,\varepsilon}^\gamma(\mathbb{R}^{n-2} 
\times \mathbb{R}^2)}\] 
and 
\[\|\mathbb{L}_\varepsilon (a,f) - (b,h)\|_{\mathcal{C}_{\mu,\varepsilon}^\gamma(\mathbb{R}^{n-
2} \times \mathbb{R}^2)} \leq C \, \Big ( \frac{\varepsilon}{\delta} + \frac{\varepsilon
^2}{\delta^2} \Big ) \, \|b\|_{\mathcal{C}_{\mu,\varepsilon}^\gamma(\mathbb{R}^{n-2} \times 
\mathbb{R}^2)}.\] 
\end{proposition}

\textbf{Proof.} 
By Corollary 3.6, there exists a pair $(a,f) \in \mathcal{E}_{\mu,\varepsilon}^{2,\gamma}
(\mathbb{R}^{n-2} \times \mathbb{R}^2)$ such that 
\[\|(a,f)\|_{\mathcal{C}_{\mu,\varepsilon}^{2,\gamma}(\mathbb{R}^{n-2} \times \mathbb{R}^2)} 
\leq C \, \varepsilon^2 \, \|(b,h)\|_{\mathcal{C}_{\mu,\varepsilon}^\gamma(\mathbb{R}^{n-2} 
\times \mathbb{R}^2)}\] 
and 
\[\mathbb{L}_\varepsilon (a,f) = (b,h).\] 
Let $\zeta$ be a cut-off function on $\mathbb{R}^{n-2}$ such that $\zeta(x) = 1$ for $|x| \leq 
\delta$, $\zeta(x) = 0$ for $|x| \geq 2\delta$, and 
\[\sup \, \delta \, |\nabla \zeta| + \sup \, \delta^2 \, |\nabla^2 \zeta| \leq C.\] 
Furthermore, let $\eta$ be a cut-off function on $\mathbb{R}^{n-2}$ satisfying $\eta(y) = 1$ for 
$|y| \leq 2\delta$, $\eta(y) = 0$ for $|y| \geq 4\delta$, and 
\[\sup \, \delta \, |\nabla \eta| + \sup \, \delta^2 \, |\nabla^2 \eta| \leq C.\] 
Then we have the estimates 
\[\|(\eta \, \zeta \, a,\eta \, \zeta \, f)\|_{\mathcal{C}_{\mu,\varepsilon}^{2,\gamma}(\mathbb
{R}^{n-2} \times \mathbb{R}^2)} \leq C \, \varepsilon^2 \, \|(b,h)\|_{\mathcal{C}_{\mu,
\varepsilon}^\gamma(\mathbb{R}^{n-2} \times \mathbb{R}^2)}\] 
and 
\begin{align*} 
&\|\mathbb{L}_\varepsilon(\eta \, \zeta \, a,\eta \, \zeta \, f) - (b,h)\|_{\mathcal{C}_{\mu,
\varepsilon}^\gamma(\mathbb{R}^{n-2} \times \mathbb{R}^2)} \\ 
&= \|\mathbb{L}_\varepsilon (\eta \, \zeta \, a,\eta \, \zeta \, f) - \eta \, \zeta \, \mathbb
{L}_\varepsilon a\|_{\mathcal{C}_{\mu,\varepsilon}^\gamma(\mathbb{R}^{n-2} \times \mathbb{R}^2)} 
\\ 
&\leq C \, \frac{1}{\varepsilon^2} \, \Big ( \frac{\varepsilon}{\delta} + \frac{\varepsilon
^2}{\delta^2} \Big ) \, \|(a,f)\|_{\mathcal{C}_{\mu,\varepsilon}^{1,\gamma}(\mathbb{R}^{n-2} 
\times \mathbb{R}^2)} \\ 
&\leq C \, \Big ( \frac{\varepsilon}{\delta} + \frac{\varepsilon^2}{\delta^2} \Big ) \, \|(b,h)
\|_{\mathcal{C}_{\mu,\varepsilon}^\gamma(\mathbb{R}^{n-2} \times \mathbb{R}^2)}. 
\end{align*} 
From this the assertion follows. \\

\section{Construction of the approximate solutions}

In this section, we define a family of approximate solutions which concentrate near $S$ as 
$\varepsilon \to 0$. To this end, we identify the total space of the normal bundle $NS$ with a 
neighborhood of the submanifold $S$ by means of the exponential map 
\[\exp: NS \to M.\] 
Note that this identification is not isometric. To see this, we choose an orthonormal basis $\{e
_i: 1 \leq i \leq n-2\}$ for the horizontal subspace, and an orthonormal basis $\{e^\perp
_\alpha: 1 \leq \alpha \leq 2\}$ for the vertical subspace. Using Jacobi's equation, we obtain 
\begin{align*} 
\exp_*(e_i) 
&= \pi_*(e_i) + \sum_{j=1}^{n-2} \sum_{\rho=1}^2 h_{ij,\rho} \, y_\rho \, \pi_*(e_j) \\ 
&- \frac{1}{2} \sum_{j=1}^{n-2} \sum_{\rho,\sigma=1}^2 R_{i\rho\sigma j} \, y_\rho \, y_\sigma 
\, \pi_*(e_j) \\ 
&- \frac{1}{2} \sum_{\beta,\rho,\sigma=1}^2 R_{i\rho\sigma\beta} \, y_\rho \, y_\sigma \, 
e^\perp_\beta + O(|y|^3). 
\end{align*} 
This implies 
\begin{align*} 
\langle \exp_*(e_i),\exp_*(e_j) \rangle 
&= \delta_{ij} + 2 \sum_{\rho=1}^2 h_{ij,\rho} \, y_\rho \\ 
&+ \sum_{k=1}^{n-2} \sum_{\rho,\sigma=1}^2 h_{ik,\rho} \, h_{jk,\sigma} \, y_\rho \, y_\sigma \\ 
&- \sum_{\rho,\sigma=1}^2 R_{i\rho\sigma j} \, y_\rho \, y_\sigma + O(|y|^3). 
\end{align*} 
Furthermore, we have 
\[\langle \exp_*(e_i),\exp_*(e^\perp_\alpha) \rangle = O(|y|^2)\] 
and 
\[\langle \exp_*(e^\perp_\alpha),\exp_*(e^\perp_\beta) \rangle = \delta_{\alpha\beta} - \frac
{1}{3} \sum_{\rho,\sigma=1}^2 R_{\alpha\rho\sigma\beta} \, y_\rho \, y_\sigma + O(|y|^3).\] 

Let $v$ be a section of the normal bundle $NS$. For every point $x \in S$, there is a unique 
one-vortex solution $(B,\psi)$ on $NS_x$ with center $v_x$. We now define a pair $(A,\phi)$ by 
\[A(e_\alpha^\perp) = B(e_\alpha^\perp),\] 
\[A(e_i) = -\nabla_i v_\rho \, B(e_\rho^\perp),\] 
and 
\[\phi = \psi.\]

\begin{proposition}
The curvature of $A$ is given by 
\[F_A(e_\alpha^\perp,e_\beta^\perp) = F_B(e_\alpha^\perp,e_\beta^\perp),\] 
\[F_A(e_i,e_\alpha^\perp) = -\nabla_i v_\rho \, F_B(e_\rho^\perp,e_\alpha^\perp),\] 
and 
\[F_A(e_i,e_j) = \nabla_i v_\rho \, \nabla_j v_\sigma \, F_A(e_\rho^\perp,e_\sigma^\perp) + 
C_{ij} + A \big ( C_{ij} \, (y - v) \big ),\] 
where $C_{ij} \in \Lambda^2 NS$ is the curvature of the normal bundle. Furthermore, the 
covariant derivative of the section $\phi$ is given by 
\[D_{A,e_\alpha^\perp} \phi = D_{B,e_\alpha^\perp} \psi\] 
and 
\[D_{A,e_i} \phi = -\nabla_i v_\rho \, D_{B,e_\rho^\perp} \psi.\] 
\end{proposition} 

\textbf{Proof.} 
Using the identity 
\[[e_i,e_j] = -C_{ij} \, y,\] 
we obtain 
\begin{align*} 
F_A(e_i,e_j) 
&= C_{ij} + \nabla_{e_i} A(e_j) - \nabla_{e_j} A(e_i) - A([e_i,e_j]) \\ 
&= C_{ij} - (\nabla_i \nabla_j v_\rho - \nabla_j \nabla_i v_\rho) \, A(e_\rho^\perp) \\ 
&+ \nabla_i v_\rho \, \nabla_j v_\sigma \, (\nabla_{e_\rho^\perp} A(e_\sigma^\perp) - \nabla_{e
_\sigma^\perp} A(e_\rho^\perp)) + A \big ( C_{ij} \, y \big ) \\ 
&= C_{ij} - A \big ( C_{ij} \, v \big ) \\ 
&+ \nabla_i v_\rho \, \nabla_j v_\sigma \, (\nabla_{e_\rho^\perp} A(e_\sigma^\perp) - \nabla_{e
_\sigma^\perp} A(e_\rho^\perp)) + A \big ( C_{ij} \, y \big ) \\ 
&= \nabla_i v_\rho \, \nabla_j v_\sigma \, F_A(e_\rho^\perp,e_\sigma^\perp) \\ 
&+ C_{ij} + A \big ( C_{ij} \, (y - v) \big ). 
\end{align*}
This proves the assertion. \\

Note that the pair $(B,\psi)$ satisfies the estimates 
\[|F_B| \leq \frac{C}{\varepsilon^2} \, e^{-\frac{\mu r}{\varepsilon}},\] 
\[|D_B \psi| \leq C \, \varepsilon^{-1} \, e^{-\frac{\mu r}{\varepsilon}},\] 
\[0 \leq 1 - |\psi|^2 \leq C \, e^{-\frac{\mu r}{\varepsilon}}\] 
for suitable constant $\mu > 0$. \\

Throughout this paper, we will assume that the normal vector field $v$ satisfies the estimate 
$\|v\|_{\mathcal{C}^{2,\gamma}(S)} \leq \varepsilon$. \\

\begin{proposition}
If the normal vector field $v$ satisfies the estimate $\|v\|_{\mathcal{C}^{2,\gamma}(S)} \leq 
\varepsilon$, then the error term verifies the estimate 
\[\Big \| \Big ( d^* F_A - \frac{1}{2\varepsilon^2} \, (\phi \, \overline{D_A \phi} - \overline
{\phi} \, D_A \phi),D_A^* D_A \phi - \frac{1}{2\varepsilon^2} \, (1 - |\phi|^2) \, \phi \Big ) 
\Big \|_{\mathcal{C}_{\mu,\varepsilon}^\gamma(M)} \leq C\] 
for some $\mu > 0$.
\end{proposition}

\textbf{Proof.} 
Since $(B,\psi)$ is a solution of the vortex equations on $\mathbb{R}^2$, we have 
\[\sum_{\rho=1}^2 \nabla_{e_\rho^\perp} F_A(e_\rho^\perp,e_\alpha^\perp) + \frac{1}{2\varepsilon
^2} \, (\phi \, \overline{D_{A,e_\alpha^\perp} \phi} - \overline{\phi} \, D_{A,e_\alpha^\perp} 
\phi) = 0\] 
and 
\[\sum_{\rho=1}^2 D_{A,e_\rho^\perp} D_{A,e_\rho^\perp} \phi + \frac{1}{2\varepsilon^2} \, (1 - 
|\phi|^2) \, \phi = 0.\] 
Furthermore, we have 
\begin{align*} 
&\sum_{\rho=1}^2 \nabla_{e_\rho^\perp} F_A(e_\rho^\perp,e_i) + \frac{1}{2\varepsilon^2} \, (\phi 
\, \overline{D_{A,e_i} \phi} - \overline{\phi} \, D_{A,e_i} \phi) \\ 
&= -\sum_{\rho=1}^2 \nabla_i v_\rho \, \Big ( \sum_{\beta=1}^2 \nabla_{e_\beta^\perp} F_A(e
_\beta^\perp,e_\rho^\perp) + \frac{1}{2\varepsilon^2} \, (\phi \, \overline{D_{A,e_\rho^\perp} 
\phi} - \overline{\phi} \, D_{A,e_\rho^\perp} \phi) \Big ) \\ 
&= 0. 
\end{align*} 
From this it follows that 
\[\Big \| \Big ( d^{*_{g_0}} F_A - \frac{1}{2\varepsilon^2} \, (\phi \, \overline{D_A 
\phi} - \overline{\phi} \, D_A \phi),D_A^{*_{g_0}} D_A \phi - \frac{1}{2\varepsilon^2} \, 
(1 - |\phi|^2) \, \phi \Big ) \Big \|_{\mathcal{C}_{\mu,\varepsilon}^\gamma(M)} \leq C.\] 
Here, $g_0$ denotes the product metric on $NS$, i.e. 
\begin{align*} 
&g_0(e_i,e_j) = \delta_{ij} \\ 
&g_0(e_i,e_\alpha^\perp) = 0 \\ 
&g_0(e_\alpha^\perp,e_\beta^\perp) = \delta_{\alpha\beta}. 
\end{align*} 
Let $g$ be the pull-back of the Riemannian metric on $M$ under the exponential map $\exp: NS \to 
M$. Then the metric $g$ satisfies an asymptotic expansion of the form 
\begin{align*} 
&g(e_i,e_j) = \delta_{ij} + 2 \sum_{\rho=1}^2 h_{ij,\rho} \, y_\rho + O(|y|^2) \\ 
&g(e_i,e_\alpha^\perp) = O(|y|^2) \\ 
&g(e_\alpha^\perp,e_\beta^\perp) = \delta_{\alpha\beta} + O(|y|^2), 
\end{align*} 
where $h$ denotes the second fundamental form of $S$. In particular, the volume form of $g$ is 
related to the volume form of $g_0$ by 
\[\bigg ( \frac{\det g}{\det g_0} \bigg )^{\frac{1}{2}} = 1 + H_\rho \, y_\rho + O(|y|^2),\] 
where $H$ is the mean curvature vector of $S$. Since the mean curvature of $S$ is $0$, we obtain 
\[\bigg ( \frac{\det g}{\det g_0} \bigg )^{\frac{1}{2}} = 1 + O(|y|^2).\] 
Thus, we conclude that 
\[\Big \| \Big ( d^* F_A - \frac{1}{2\varepsilon^2} \, (\phi \, \overline{D_A \phi} - 
\overline{\phi} \, D_A \phi),D_A^* D_A \phi - \frac{1}{2\varepsilon^2} \, (1 - |\phi|^2) 
\, \phi \Big ) \Big \|_{\mathcal{C}_{\mu,\varepsilon}^\gamma(M)} \leq C.\] 
This proves the assertion. \\

\section{Estimates for the operator $\mathbb{L}_\varepsilon = L_\varepsilon + T_\varepsilon 
T_\varepsilon^*$ in weighted H\"older spaces}

Our aim in this section is to analyze the mapping properties of the linearized operator 
\begin{align*} 
\mathbb{L}_\varepsilon: \; &\Omega^1(\mathbb{R}^{n-2} \times \mathbb{R}^2,i\mathbb{R}) \oplus 
\Omega^0(\mathbb{R}^{n-2} \times \mathbb{R}^2,L) \\ 
&\to \Omega^1(\mathbb{R}^{n-2} \times \mathbb{R}^2,i\mathbb{R}) \oplus \Omega^0(\mathbb{R}^{n-2} 
\times \mathbb{R}^2,L). 
\end{align*}

\begin{proposition}
Suppose that $(b,h) \in \mathcal{C}_{\mu,\varepsilon}^\gamma(M)$ is supported in the set $\{p 
\in M: \dist(p,S) \leq 2\delta\}$ and satisfies 
\[\int_{NS_x} \varepsilon^2 \, \sum_{\alpha=1}^4 \langle b(e_\alpha^\perp),F_A(w,e_\alpha^\perp) 
\rangle + \int_{NS_x} \langle h,D_{A,w} \phi \rangle = 0\] 
for all $x \in S$ and all $w \in NS$. Then there exists a pair $(a,f) \in \mathcal{C}_{\mu,
\varepsilon}^{2,\gamma}(M)$ which is supported in the region $\{p \in M: \dist(p,S) \leq 4\delta
\}$ such that 
\[\|(a,f)\|_{\mathcal{C}_{\mu,\varepsilon}^{2,\gamma}(M)} \leq C \, \varepsilon^2 \, \|b\|
_{\mathcal{C}_{\mu,\varepsilon}^\gamma(M)}\] 
and 
\[\|\mathbb{L}_\varepsilon (a,f) - (b,h)\|_{\mathcal{C}_{\mu,\varepsilon}^\gamma(M)} \leq C \, 
\Big ( \delta + \frac{\varepsilon}{\delta} + \frac{\varepsilon^2}{\delta^2} \Big ) \, \|(b,h)\|
_{\mathcal{C}_{\mu,\varepsilon}^\gamma(M)}.\]
\end{proposition}

\textbf{Proof.} 
Let $\{\zeta^{(j)}: 1 \leq j \leq j_0\}$ be a partition of unity on $S$ such that each function 
$\zeta^{(j)}$ is supported in a ball $B_\delta(p_j)$, and 
\[|\{1 \leq j \leq j_0: x \in B_{4\delta}(p_j)\}| \leq C\] 
for all $x \in S$ and some uniform constant $C$. For each $1 \leq j \leq j_0$, there exists a 
pair $(a^{(j)},f^{(j)}) \in \mathcal{C}_{\mu,\varepsilon}^{2,\gamma}(M)$ which is supported in 
the region $\{(x,y) \in NS: x \in B_{2\delta}(p_j), \, |y| \leq 4\delta\}$ such that 
\[\|(a^{(j)},f^{(j)})\|_{\mathcal{C}_{\mu,\varepsilon}^{2,\gamma}(M)} \leq C \, \varepsilon^2 \, 
\|(\zeta^{(j)} \, b,\zeta^{(j)} \, h)\|_{\mathcal{C}_{\mu,\varepsilon}^\gamma(M)}\] 
and 
\[\|\mathbb{L}_\varepsilon (a^{(j)},f^{(j)}) - (\zeta^{(j)} \, b,\zeta^{(j)} \, h)\|_{\mathcal
{C}_{\mu,\varepsilon}^\gamma(M)} \leq C \, \Big ( \delta + \frac{\varepsilon}{\delta} + \frac
{\varepsilon^2}{\delta^2} \Big ) \, \|(\zeta^{(j)} \, b,\zeta^{(j)} \, h)\|_{\mathcal{C}_{\mu,
\varepsilon}^\gamma(M)}.\] 
We now define 
\[(a,f) = \sum_{j=1}^{j_0} (a^{(j)},f^{(j)}).\] 
Then we have the estimates 
\begin{align*} 
\|(a,f)\|_{\mathcal{C}_{\mu,\varepsilon}^{2,\gamma}(M)} 
&\leq C \, \sup_{1 \leq j \leq j_0} \|(a^{(j)},f^{(j)})\|_{\mathcal{C}_{\mu,\varepsilon}^{2,
\gamma}(M)} \\ 
&\leq C \, \varepsilon^2 \, \sup_{1 \leq j \leq j_0} \|(\zeta^{(j)} \, b,\zeta^{(j)} \, h)\|
_{\mathcal{C}_{\mu,\varepsilon}^\gamma(M)} \\ 
&\leq C \, \varepsilon^2 \, \|(b,h)\|_{\mathcal{C}_{\mu,\varepsilon}^\gamma(M)} 
\end{align*} 
and 
\begin{align*} 
\|\mathbb{L}_\varepsilon (a,f) - (b,h)\|_{\mathcal{C}_{\mu,\varepsilon}^\gamma(M)} 
&\leq C \, \sup_{1 \leq j \leq j_0} \|\mathbb{L}_\varepsilon (a^{(j)},f^{(j)}) - (\zeta^{(j)} \, 
b,\zeta^{(j)} \, h)\|_{\mathcal{C}_{\mu,\varepsilon}^\gamma(M)} \\ 
&\leq C \, \Big ( \delta + \frac{\varepsilon}{\delta} + \frac{\varepsilon^2}{\delta^2} \Big ) \, 
\sup_{1 \leq j \leq j_0} \|(\zeta^{(j)} \, b,\zeta^{(j)} \, h)\|_{\mathcal{C}_{\mu,\varepsilon}
^\gamma(M)} \\ 
&\leq C \, \Big ( \delta + \frac{\varepsilon}{\delta} + \frac{\varepsilon^2}{\delta^2} \Big ) \, 
\|(b,h)\|_{\mathcal{C}_{\mu,\varepsilon}^\gamma(M)}. 
\end{align*} 
This proves the assertion. \\

\begin{proposition}
For every $(b,h) \in \mathcal{C}_{\mu,\varepsilon}^\gamma(M)$, there exists a pair $(a,f) \in 
\mathcal{C}_{\mu,\varepsilon}^{2,\gamma}(M)$ such that 
\[\|(a,f)\|_{\mathcal{C}_{\mu,\varepsilon}^{2,\gamma}(M)} \leq C \, \varepsilon^2 \, \|(b,h)\|
_{\mathcal{C}_{\mu,\varepsilon}^\gamma(M)}\] 
and 
\[\Big ( \nabla^* \nabla a + \frac{1}{\varepsilon^2} \, a,\nabla^* \nabla f + \frac
{1}{\varepsilon^2} \, f \Big ) = (b,h).\] 
\end{proposition}

\textbf{Proof.} 
By Schauder estimates, it suffices to show that 
\[\sup \, e^{\frac{\mu \, \dist(p,S)}{\varepsilon}} \, (\varepsilon \, |a| + |f|) \leq C \, \sup 
\, \varepsilon^2 \, e^{\frac{\mu \, \dist(p,S)}{\varepsilon}} \, 
\Big ( \varepsilon \, \big | \nabla^* \nabla a + \frac{1}{\varepsilon^2} \, a \big | 
+ \big | \nabla^* \nabla f + \frac{1}{\varepsilon^2} \, f \big | \Big ).\] 
Suppose that there exists a sequence of positive real numbers $\varepsilon_j$ and a sequence of 
pairs $(a^{(j)},f^{(j)}) \in \mathcal{C}_{\mu,\varepsilon}^{2,\gamma}(M)$ such that 
\[\sup \, e^{\frac{\mu \, \dist(p,S)}{\varepsilon_j}} \, (\varepsilon_j \, |a^{(j)}| + |f^{(j)}|) 
= 1\] 
and 
\[\sup \, \varepsilon_j^2 \, e^{\frac{\mu \, \dist(p,S)}{\varepsilon_j}} \, \Big ( \varepsilon_j 
\, \big | \nabla^* \nabla a^{(j)} + \frac{1}{\varepsilon_j^2} \, a^{(j)} \big | + \big | \nabla
^* \nabla f^{(j)} + \frac{1}{\varepsilon_j^2} \, f^{(j)} \big | \Big ) \to 0.\] 
Then there exists a sequence of points $p_j \in M$ such that 
\[e^{\frac{\mu \, \dist(p_j,S)}{\varepsilon_j}} \, (\varepsilon_j \, |a^{(j)}(p_j)| + |f^{(j)}(p
_j)|) \geq \frac{1}{2}.\] 
After rescaling, we obtain a sequence of pairs $(\tilde{a}^{(j)},\tilde{f}^{(j)})$ such that 
\[\sup \, e^{-\mu \, \dist(p,p_j)} \, (|\tilde{a}^{(j)}| + |\tilde{f}^{(j)}|) \leq 1\] 
and 
\[\sup \, e^{-\mu \, \dist(p,p_j)} \, (|\nabla^* \nabla \tilde{a}^{(j)} + \tilde{a}^{(j)}| + 
|\nabla^* \nabla \tilde{f}^{(j)} + \tilde{f}^{(j)}|) \to 0.\] 
Moreover, we have 
\[|\tilde{a}^{(j)}(p_j)| + |\tilde{f}^{(j)}(p_j)| \geq \frac{1}{2}.\] 
Taking the limit as $j \to \infty$, we obtain a pair $(\tilde{a},\tilde{f})$ such that 
\[\sup \, e^{-\mu \, \dist(p,p_0)} \, (|\tilde{a}| + |\tilde{f}|) \leq 1\] 
and 
\[(\nabla^* \nabla \tilde{a} + \tilde{a},\nabla^* \nabla \tilde{f} + \tilde{f}) = 0.\] 
If $\mu$ is sufficiently small, we conclude that $(\tilde{a},\tilde{f}) = 0$. This is a 
contradiction. \\

\begin{proposition}
Suppose that $(b,h) \in \mathcal{C}_{\mu,\varepsilon}^\gamma(M)$ is supported in the region $\{p 
\in M: \dist(p,S) \geq \delta\}$. Then there exists a pair $(a,f) \in \mathcal{C}_{\mu,
\varepsilon}^{2,\gamma}(M)$ which is supported in the region $\{p \in M: \dist(p,S) \geq \frac
{\delta}{2}\}$ such that 
\[\|(a,f)\|_{\mathcal{C}_{\mu,\varepsilon}^{2,\gamma}(M)} \leq C \, \varepsilon^2 \, \|(b,h)\|
_{\mathcal{C}_{\mu,\varepsilon}^\gamma(M)}\] 
and 
\[\|\mathbb{L}_\varepsilon (a,f) - (b,h)\|_{\mathcal{C}_{\mu,\varepsilon}^\gamma(M)} \leq C \, 
\Big ( \frac{\varepsilon}{\delta} + \frac{\varepsilon^2}{\delta^2} + e^{-\frac{\mu 
\delta}{\varepsilon}} \Big ) \, \|(b,h)\|_{\mathcal{C}_{\mu,\varepsilon}^\gamma(M)}.\]
\end{proposition}

\textbf{Proof.} 
By Proposition 5.2, we can find a pair $(a,f)$ such that 
\[\|(a,f)\|_{\mathcal{C}_{\mu,\varepsilon}^{2,\gamma}(M)} \leq C \, \varepsilon^2 \, \|(b,h)\|
_{\mathcal{C}_{\mu,\varepsilon}^\gamma(M)}\] 
and 
\[\Big ( \nabla^* \nabla a + \frac{1}{\varepsilon^2} \, a,\nabla^* \nabla f + \frac
{1}{\varepsilon^2} \, f \Big ) = (b,h).\] 
Let $\eta$ be a cut-off function such that $\eta(p) = 0$ for $\dist(p,S) \leq \frac
{\delta}{2}$, $\eta(p) = 1$ for $\dist(p,S) \geq \delta$ and 
\[\sup \, \delta \, |\nabla \eta| + \sup \, \delta^2 \, |\nabla^2 \eta| \leq C.\] 
Then the pair $(\eta \, a,\eta \, f)$ is supported in the region $\{p \in M: \dist(p,S) \geq 
\frac{\delta}{2}\}$ and satisfies 
\begin{align*} 
&\|\mathbb{L}_\varepsilon (\eta \, a,\eta \, f) - (b,h)\|_{\mathcal{C}_{\mu,\varepsilon}^\gamma
(M)} \\ 
&\leq \Big \| \mathbb{L}_\varepsilon (\eta \, a,\eta \, f) - \Big ( \nabla^* \nabla (\eta \, a) 
+ \frac{1}{\varepsilon^2} \, \eta \, a,\nabla^* \nabla (\eta \, f) + \frac{1}{\varepsilon^2} \, 
\eta \, f \Big ) \Big \|_{\mathcal{C}_{\mu,\varepsilon}^\gamma(M)} \\ 
&+ \Big \| (\nabla^* \nabla (\eta \, a) - \eta \, \nabla^* \nabla a,\nabla^* \nabla (\eta \, f) 
- \eta \, \nabla^* \nabla f) \Big \|_{\mathcal{C}_{\mu,\varepsilon}^\gamma(M)} \\ 
&\leq C \, \frac{1}{\varepsilon^2} \, \Big ( e^{-\frac{\mu \delta}{\varepsilon}} + \frac
{\varepsilon}{\delta} + \frac{\varepsilon^2}{\delta^2} \Big ) \, \|(a,f)\|_{\mathcal{C}_{\mu,
\varepsilon}^{2,\gamma}(M)} \\ 
&\leq C \, \Big ( e^{-\frac{\mu \delta}{\varepsilon}} + \frac{\varepsilon}{\delta} + \frac
{\varepsilon^2}{\delta^2} \Big ) \, \|(b,h)\|_{\mathcal{C}_{\mu,\varepsilon}^{2,\gamma}(M)}. 
\end{align*} 
This proves the assertion. \\

In the following, we will choose $\delta = \varepsilon^{\frac{1}{2}}$. Let $\kappa$ be a cut-off 
function such that $\kappa(p) = 1$ for $\dist(p,S) \leq \varepsilon^{\frac{1}{2}}$ and $\kappa
(p) = 0$ for $\dist(p,S) \geq 2\varepsilon^{\frac{1}{2}}$. \\

Let $\mathcal{E}_{\mu,\varepsilon}^{k,\gamma}(M)$ be the set of all pairs $(b,h) \in \Omega^1(M,
i\mathbb{R}) \oplus \Omega^0(M,L)$ such that $(b,h) \in \mathcal{C}_{\mu,\varepsilon}^{k,\gamma}
(M)$ and 
\[\int_{NS_x} \varepsilon^2 \, \kappa \, \sum_{\alpha=1}^4 \langle b(e_\alpha^\perp),F_A(w,e
_\alpha^\perp) \rangle + \int_{NS_x} \kappa \, \langle h,D_{A,w} \phi \rangle = 0\] 
for all $x \in S$ and $w \in NS_x$. \\

We denote by $\Id - \mathbb{P}$ the fibrewise projection from $\mathcal{C}_\nu^\gamma(M)$ to the 
subspace $\mathcal{E}_\nu^\gamma(M)$. Hence, for every pair $(b,h)$ there exists a normal vector 
field $w$ such that 
\[\mathbb{P} (b,h) = (F_A(w,\cdot),D_{A,w} \phi).\] 
Let $\Pi$ be the linear operator which assigns to every pair $(b,h)$ the vector field 
\[\Pi (b,h) = w.\] 
It is not difficult to show that 
\[\|\Pi (b,h)\|_{\mathcal{C}^\gamma(S)} \leq C \, \varepsilon^{1-\gamma} \, \|(b,h)\|_{\mathcal
{C}_{\mu,\varepsilon}^\gamma(M)}\] 
and 
\[\|\mathbb{P} (b,h)\|_{\mathcal{C}_{\mu,\varepsilon}^\gamma(M)} \leq C \, \|(b,h)\|_{\mathcal
{C}_{\mu,\varepsilon}^\gamma(M)}.\]

\begin{proposition}
For every pair $(b,h) \in \mathcal{E}_{\mu,\varepsilon}^\gamma(M)$ there exists a pair $(a,f) 
\in \mathcal{C}_{\mu,\varepsilon}^{2,\gamma}(M)$ such that 
\[\|(a,f)\|_{\mathcal{C}_{\mu,\varepsilon}^{2,\gamma}(M)} \leq C \, \varepsilon^2 \, \|(b,h)\|
_{\mathcal{C}_{\mu,\varepsilon}^\gamma(M)}\] 
and 
\[\|\mathbb{L}_\varepsilon (a,f) - (b,h)\|_{\mathcal{C}_{\mu,\varepsilon}^\gamma(M)} \leq C \, 
\varepsilon^{\frac{1}{2}} \, \|(b,h)\|_{\mathcal{C}_{\mu,\varepsilon}^\gamma(M)}.\]
\end{proposition}

\textbf{Proof.} 
Apply Proposition 5.1 to $(\kappa \, b,\kappa \, h)$ and Proposition 5.3 to $((1 - \kappa) \, b,
(1 - \kappa) \, h)$. \\

\begin{proposition}
For every $(b,h) \in \mathcal{E}_{\mu,\varepsilon}^\gamma(M)$ there exists a pair $(a,f) \in 
\mathcal{C}_{\mu,\varepsilon}^{2,\gamma}(M)$ such that 
\[\|(a,f)\|_{\mathcal{C}_{\mu,\varepsilon}^{2,\gamma}(M)} \leq C \, \varepsilon^2 \, \|(b,h)\|
_{\mathcal{C}_{\mu,\varepsilon}^\gamma(M)}\] 
and 
\[(\Id - \mathbb{P}) \, \mathbb{L}_\varepsilon (a,f) = (b,h).\] 
Furthermore, the pair $(a,f)$ satisfies the estimate 
\[\|\Pi \, \mathbb{L}_\varepsilon (a,f)\|_{\mathcal{C}^\gamma(S)} \leq C \, \varepsilon^{\frac
{5}{4}} \, \|(b,h)\|_{\mathcal{C}_{\mu,\varepsilon}^\gamma(M)}.\] 
\end{proposition}

\textbf{Proof.} 
By Proposition 5.4, there exists an operator $\mathbb{S}: \mathcal{E}_{\mu,\varepsilon}^\gamma
(M) \to \mathcal{C}_{\mu,\varepsilon}^{2,\gamma}(M)$ such that 
\[\|\mathbb{S} (b,h)\|_{\mathcal{C}_{\mu,\varepsilon}^{2,\gamma}(M)} \leq C \, \varepsilon^2 \, 
\|(b,h)\|_{\mathcal{C}_{\mu,\varepsilon}^\gamma(M)}\] 
and 
\[\|\mathbb{L}_\varepsilon \, \mathbb{S} (b,h) - (b,h)\|_{\mathcal{C}_{\mu,\varepsilon}^\gamma
(M)} \leq C \, \varepsilon^{\frac{1}{2}} \, \|(b,h)\|_{\mathcal{C}_{\mu,\varepsilon}^\gamma(M)}.
\] 
This implies 
\[\|\Pi \, \mathbb{L}_\varepsilon \, \mathbb{S} (b,h)\|_{\mathcal{C}^\gamma(S)} = \|\Pi(\mathbb
{L}_\varepsilon \, \mathbb{S} (b,h) - (b,h))\|_{\mathcal{C}^\gamma(S)} \leq C \, \varepsilon
^{\frac{3}{2}-\gamma} \, \|(b,h)\|_{\mathcal{C}_{\mu,\varepsilon}^\gamma(M)}.\] 
From this it follows that 
\[\|(\Id - \mathbb{P}) \, \mathbb{L}_\varepsilon \, \mathbb{S} (b,h) - (b,h)\|_{\mathcal{C}
_{\mu,\varepsilon}^\gamma(M)} \leq C \, \varepsilon^{\frac{1}{2}} \, \|(b,h)\|_{\mathcal{C}
_{\mu,\varepsilon}^\gamma(M)}.\] 
Therefore, the operator $(\Id - \mathbb{P}) \, \mathbb{L}_\varepsilon \, \mathbb{S}: \mathcal{E}
_{\mu,\varepsilon}^\gamma(M) \to \mathcal{E}_{\mu,\varepsilon}^\gamma(M)$ is invertible. Hence, 
if we define 
\[(a,f) = \mathbb{S} \, \big [ (\Id - \mathbb{P}) \, \mathbb{L}_\varepsilon \, \mathbb{S} \big ]
^{-1} \, (b,h),\] 
then $(a,f)$ satisfies 
\[\|(a,f)\|_{\mathcal{C}_{\mu,\varepsilon}^{2,\gamma}(M)} \leq C \, \varepsilon^2 \, \|b\|
_{\mathcal{C}_{\mu,\varepsilon}^\gamma(M)}\] 
and 
\[(\Id - \mathbb{P}) \, \mathbb{L}_\varepsilon (a,f) = (b,h).\] 
This proves the assertion. \\

\section{From approximate to exact solutions}

\begin{proposition}
Let $v$ be a normal vector field along $S$, and let $A$ be the approximate solution associated 
to $v$. Then there exists a pair $(\tilde{A},\tilde{\phi}) = (A + a,\phi + f)$ such that $(a,f) 
\in \mathcal{C}_{\mu,\varepsilon}^{2,\gamma}(M)$, 
\[\|(a,f)\|_{\mathcal{C}_{\mu,\varepsilon}^{2,\gamma}(M)} \leq C \, \varepsilon^2\] 
and 
\[(\Id - \mathbb{P}) \Big ( d^* F_{\tilde{A}} - \frac{1}{2\varepsilon^2} \, (\tilde{\phi} \, 
\overline{D_{\tilde{A}} \tilde{\phi}} - \overline{\tilde{\phi}} \, D_{\tilde{A}} \tilde{\phi}) + 
\frac{1}{\varepsilon} \, du,D_{\tilde{A}}^* D_{\tilde{A}} \tilde{\phi} - \frac{1}{2\varepsilon
^2} \, (1 - |\tilde{\phi}|^2) \, \tilde{\phi} - \frac{1}{\varepsilon} \, \tilde{\phi} \, u \Big 
) = 0,\] 
where 
\[u = \varepsilon \, d^* a + \frac{1}{2\varepsilon} \, (\phi \, \overline{f} - \overline{\phi} 
\, f).\] 
Furthermore, the pair $(a,f)$ satisfies the estimate 
\[\|\Pi \, \mathbb{L}_\varepsilon (a,f)\|_{\mathcal{C}^\gamma(S)} \leq C \, \varepsilon^{\frac
{5}{4}}.\] 
\end{proposition}

\textbf{Proof.}
We use the identity 
\begin{align*} 
&\Big ( d^* F_{\tilde{A}} - \frac{1}{2\varepsilon^2} \, (\tilde{\phi} \, \overline{D_{\tilde{A}} 
\tilde{\phi}} - \overline{\tilde{\phi}} \, D_{\tilde{A}} \tilde{\phi}),D_{\tilde{A}}^* D_{\tilde
{A}} \tilde{\phi} - \frac{1}{2\varepsilon^2} \, (1 - |\tilde{\phi}|^2) \, \tilde{\phi} - \frac
{1}{\varepsilon} \, \tilde{\phi} \, u \Big ) \\ &= \Big ( d^* F_A - \frac{1}{2\varepsilon^2} \, 
(\phi \, \overline{D_A \phi} - \overline{\phi} \, D_A \phi) + \frac{1}{\varepsilon} \, du,D_A^* 
D_A \phi - \frac{1}{2\varepsilon^2} \, (1 - |\phi|^2) \, \phi \Big ) \\ 
&+ L_\varepsilon (a,f) + Q(a,f). 
\end{align*} 
Here, the remainder term $Q(a,f)$ is at least quadratic in $(a,f)$. This implies 
\[\|Q(a,f)\|_{\mathcal{C}_{\mu,\varepsilon}^\gamma(M)} \leq C \, \varepsilon^{-2} \, \|(a,f)\|
_{\mathcal{C}_{\mu,\varepsilon}^{2,\gamma}(M)}^2.\] 
Hence, if we define 
\[u = T_\varepsilon^* (a,f),\] 
then we obtain 
\begin{align*} 
&\Big ( d^* F_{\tilde{A}} - \frac{1}{2\varepsilon^2} \, (\tilde{\phi} \, \overline{D_{\tilde{A}} 
\tilde{\phi}} - \overline{\tilde{\phi}} \, D_{\tilde{A}} \tilde{\phi}) + \frac{1}{\varepsilon} 
\, du,D_{\tilde{A}}^* D_{\tilde{A}} \tilde{\phi} - \frac{1}{2\varepsilon^2} \, (1 - |\tilde
{\phi}|^2) \, \tilde{\phi} - \frac{1}{\varepsilon} \, \tilde{\phi} \, u \Big ) \\ 
&= \Big ( d^* F_A - \frac{1}{2\varepsilon^2} \, (\phi \, \overline{D_A \phi} - \overline{\phi} 
\, D_A \phi),D_A^* D_A \phi - \frac{1}{2\varepsilon^2} \, (1 - |\phi|^2) \, \phi \Big ) \\ 
&+ \mathbb{L}_\varepsilon (a,f) + Q(a,f). 
\end{align*} 
According to Proposition 5.5, there exists an operator $\mathbb{G}: \mathcal{E}_{\mu,
\varepsilon}^\gamma(M) \to \mathcal{C}_{\mu,\varepsilon}^{2,\gamma}(M)$ such that 
\[\|\mathbb{G} (b,h)\|_{\mathcal{C}_{\mu,\varepsilon}^{2,\gamma}(M)} \leq C \, \varepsilon^2 \, 
\|(b,h)\|_{\mathcal{C}_{\mu,\varepsilon}^\gamma(M)}\] 
and 
\[(\Id - \mathbb{P}) \, \mathbb{L}_\varepsilon \, \mathbb{G} = \Id.\] 
We now define a mapping $\Phi: \mathcal{C}_{\mu,\varepsilon}^{2,\gamma}(M) \to \mathcal{C}_{\mu,
\varepsilon}^{2,\gamma}(M)$ by 
\begin{align*} 
&\Phi(a,f) \\ 
&= -\mathbb{G} \, (\Id - \mathbb{P}) \,  \Big ( d^* F_A - \frac{1}{2\varepsilon^2} \, (\phi \, 
\overline{D_A \phi} - \overline{\phi} \, D_A \phi),D_A^* D_A \phi - \frac{1}{2\varepsilon^2} \, 
(1 - |\phi|^2) \, \phi \Big ) \\ &- \mathbb{G} \, (\Id - \mathbb{P}) \, Q(a,f). 
\end{align*} 
Then we have the estimate 
\begin{align*} 
&\|\Phi(a,f)\|_{\mathcal{C}_{\mu,\varepsilon}^{2,\gamma}(M)} \\ 
&\leq C \, \varepsilon^2 \, \Big \| (\Id - \mathbb{P}) \, \Big ( d^* F_A - \frac{1}{2\varepsilon
^2} \, (\phi \, \overline{D_A \phi} - \overline{\phi} \, D_A \phi),D_A^* D_A \phi - \frac{1}{2
\varepsilon^2} \, (1 - |\phi|^2) \, \phi \Big ) \Big \|_{\mathcal{C}_{\mu,\varepsilon}^\gamma
(M)} \\ 
&+ C \, \varepsilon^2 \, \|(\Id - \mathbb{P}) \, Q(a,f)\|_{\mathcal{C}_{\mu,\varepsilon}^\gamma
(M)} \\ 
&\leq C \, \varepsilon^2 \, \Big \| \Big ( d^* F_A - \frac{1}{2\varepsilon^2} \, (\phi \, 
\overline{D_A \phi} - \overline{\phi} \, D_A \phi),D_A^* D_A \phi - \frac{1}{2\varepsilon^2} \, 
(1 - |\phi|^2) \, \phi \Big ) \Big \|_{\mathcal{C}_{\mu,\varepsilon}^\gamma(M)} \\ 
&+ C \, \varepsilon^2 \, \|Q(a,f)\|_{\mathcal{C}_{\mu,\varepsilon}^\gamma(M)} \\ 
&\leq C \, \varepsilon^2 
\end{align*} 
for all $(a,f) \in \mathcal{C}_{\mu,\varepsilon}^{2,\gamma}(M)$ satisfying 
\[\|a\|_{\mathcal{C}_{\mu,\varepsilon}^{2,\gamma}(M)} \leq \varepsilon^{\frac{3}{2}}.\] 
Moreover, we have 
\begin{align*} 
\|\Phi(a,f) - \Phi(a',f')\|_{\mathcal{C}_{\mu,\varepsilon}^{2,\gamma}(M)} 
&\leq C \, \varepsilon^2 \, \|Q(a,f) - Q(a',f')\|_{\mathcal{C}_{\mu,\varepsilon}^\gamma(M)} \\ 
&\leq C \, \varepsilon^{\frac{3}{2}} \, \|(a,f) - (a',f')\|_{\mathcal{C}_{\mu,\varepsilon}^{2,
\gamma}(M)} 
\end{align*} 
for all $(a,f),(a',f') \in \mathcal{C}_{\mu,\varepsilon}^{2,\gamma}(M)$ satisfying 
\[\|(a,f)\|_{\mathcal{C}_{\mu,\varepsilon}^{2,\gamma}(M)}, \, \|(a',f')\|_{\mathcal{C}_{\mu,
\varepsilon}^{2,\gamma}(M)} \leq \varepsilon^{\frac{3}{2}}.\] 
Hence, it follows from the contraction mapping principle that there exists a pair $(a,f) \in 
\mathcal{C}_{\mu,\varepsilon}^{2,\gamma}(M)$ such that 
\[\|(a,f)\|_{\mathcal{C}_{\mu,\varepsilon}^{2,\gamma}(M)} \leq C\] 
and 
\[\Phi(a,f) = (a,f).\] 
From this it follows that 
\begin{align*} 
&\mathbb{G} \, (\Id - \mathbb{P}) \,  \Big ( d^* F_A - \frac{1}{2\varepsilon^2} \, (\phi \, 
\overline{D_A \phi} - \overline{\phi} \, D_A \phi),D_A^* D_A \phi - \frac{1}{2\varepsilon^2} \, 
(1 - |\phi|^2) \, \phi \Big ) \\ &+ (a,f) + \mathbb{G} \, (\Id - \mathbb{P}) \, Q(a,f) = 0, 
\end{align*} 
hence 
\begin{align*} 
&(\Id - \mathbb{P}) \,  \Big ( d^* F_A - \frac{1}{2\varepsilon^2} \, (\phi \, \overline{D_A 
\phi} - \overline{\phi} \, D_A \phi),D_A^* D_A \phi - \frac{1}{2\varepsilon^2} \, (1 - |\phi|^2) 
\, \phi \Big ) \\ 
&+ (\Id - \mathbb{P}) \, \mathbb{L}_\varepsilon (a,f) + (\Id - \mathbb{P}) \, Q(a,f) = 0. 
\end{align*} 
Thus, we conclude that 
\begin{align*} 
&(\Id - \mathbb{P}) \, \Big ( d^* F_{\tilde{A}} - \frac{1}{2\varepsilon^2} \, (\tilde{\phi} \, 
\overline{D_{\tilde{A}} \tilde{\phi}} - \overline{\tilde{\phi}} \, D_{\tilde{A}} \tilde{\phi}) + 
\frac{1}{\varepsilon} \, du,D_{\tilde{A}}^* D_{\tilde{A}} \tilde{\phi} - \frac{1}{2\varepsilon
^2} \, (1 - |\tilde{\phi}|^2) \, \tilde{\phi} - \frac{1}{\varepsilon} \, \tilde{\phi} \, u \Big 
) \\ 
&= 0.
\end{align*} 
This proves the assertion. \\

\begin{proposition}
If 
\[\mathbb{P} \Big ( d^* F_{\tilde{A}} - \frac{1}{2\varepsilon^2} \, (\tilde{\phi} \, \overline{D
_{\tilde{A}} \tilde{\phi}} - \overline{\tilde{\phi}} \, D_{\tilde{A}} \tilde{\phi}) + \frac
{1}{\varepsilon} \, du,D_{\tilde{A}}^* D_{\tilde{A}} \tilde{\phi} - \frac{1}{2\varepsilon^2} \, 
(1 - |\tilde{\phi}|^2) \, \tilde{\phi} - \frac{1}{\varepsilon} \, \tilde{\phi} \, u \Big ) = 0,
\] 
then $(\tilde{A},\tilde{\phi})$ is a solution of the Ginzburg-Landau equations.
\end{proposition}

\textbf{Proof.}
It follows from the definition of the pair $(\tilde{A},\tilde{\phi})$ that 
\[d^* F_{\tilde{A}} = \frac{1}{2\varepsilon^2} \, (\tilde{\phi} \, \overline{D_{\tilde{A}} 
\tilde{\phi}} - \overline{\tilde{\phi}} \, D_{\tilde{A}} \tilde{\phi}) - \frac{1}{\varepsilon} 
\, du\] 
and 
\[D_{\tilde{A}}^* D_{\tilde{A}} \tilde{\phi} = \frac{1}{2\varepsilon^2} \, (1 - |\tilde{\phi}|
^2) \, \tilde{\phi} + \frac{1}{\varepsilon} \, \tilde{\phi} \, u,\] 
where $u$ satisfies $\overline{u} = -u$. From this it follows that 
\begin{align*} 
0 
&= -\varepsilon \, d^* d^* F_{\tilde{A}} \\ 
&= -\frac{1}{2\varepsilon} \, (\tilde{\phi} \, \overline{D_{\tilde{A}}^* D_{\tilde{A}} \tilde
{\phi}} - \overline{\tilde{\phi}} \, D_{\tilde{A}}^* D_{\tilde{A}} \tilde{\phi}) + d^* du \\ 
&= -\frac{1}{2\varepsilon^2} \, (\tilde{\phi} \, \overline{\tilde{\phi} \, u} - \overline{\tilde
{\phi}} \, \tilde{\phi} \, u) + d^* du \\ 
&= \frac{1}{\varepsilon^2} \, |\tilde{\phi}|^2 \, u + d^* du. 
\end{align*} 
Thus, we conclude that $u = 0$. Hence, $(\tilde{A},\tilde{\phi})$ is a solution of the 
Ginzburg-Landau equations. \\

\section{The balancing condition}

\begin{proposition}
Let $g_0$ be the product metric on the normal bundle $NS$ (cf. Section 4). Then we have the 
identity 
\begin{align*} 
&\Pi \Big ( d^{*_{g_0}} F_A - \frac{1}{2\varepsilon^2} \, (\phi \, \overline{D_A \phi} - 
\overline{\phi} \, D_A \phi),D_A^{*_{g_0}} D_A \phi - \frac{1}{2\varepsilon^2} \, (1 - |\phi|^2) 
\, \phi \Big ) \\ 
&= \Delta v. 
\end{align*}
\end{proposition}

\textbf{Proof.} 
Using the results from Section 4, we obtain 
\begin{align*} 
&\sum_{\beta=1}^2 \nabla_{e_\beta^\perp} F_A(e_\beta^\perp,e_i) + \frac{1}{2\varepsilon^2} \, 
(\phi \, \overline{D_{A,e_i} \phi} - \overline{\phi} \, D_{A,e_i} \phi) \\ 
&= -\sum_{\rho=1}^2 \nabla_i v_\rho \, \Big ( \sum_{\beta=1}^2 \nabla_{e_\beta^\perp} F_A(e
_\beta^\perp,e_\rho^\perp) + \frac{1}{2\varepsilon^2} \, (\phi \, \overline{D_{A,e_\rho^\perp} 
\phi} - \overline{\phi} \, D_{A,e_\rho^\perp} \phi) \Big ) \\ 
&= 0. 
\end{align*} 
The Bianchi identity implies that 
\[\nabla_{e_\beta^\perp} F_A(e_i,e_\alpha^\perp) - \nabla_{e_\alpha^\perp} F_A(e_i,e_\beta
^\perp) + \nabla_{e_i} F_A(e_\alpha^\perp,e_\beta^\perp) = 0.\] 
Furthermore, we have 
\[D_{A,e_i} D_{A,e_\alpha^\perp} \phi - D_{A,e_\alpha^\perp} D_{A,e_i} \phi = F_A(e_i,e_\alpha
^\perp) \, \phi.\] 
From this it follows that 
\begin{align*} 
&\sum_{i=1}^{n-2} \sum_{\alpha,\beta=1}^2 \varepsilon^2 \, \big ( \nabla_{e_\beta^\perp} \langle 
F_A(e_i,e_\beta^\perp),F_A(e_i,e_\alpha^\perp) \rangle - \frac{1}{2} \, \nabla_{e_\alpha^\perp} 
\langle F_A(e_i,e_\beta^\perp),F_A(e_i,e_\beta^\perp) \rangle \\ 
&+ \nabla_{e_i} \langle F_A(e_i,e_\beta^\perp),F_A(e_\alpha^\perp,e_\beta^\perp) \rangle \big ) 
\, w^\alpha \\ 
&+ \sum_{i=1}^{n-2} \sum_{\alpha=1}^2 \big ( -\frac{1}{2} \, \nabla_{e_\alpha^\perp} \langle D
_{A,e_i} \phi,D_{A,e_i} \phi \rangle + \nabla_{e_i} \langle D_{A,e_i} \phi,D_{A,e_\alpha^\perp} 
\phi \rangle \big ) \, w^\alpha \\ 
&= \sum_{i=1}^{n-2} \sum_{\alpha,\beta=1}^2 \varepsilon^2 \, \langle D_{A,e_i} F_A(e_i,e_\beta
^\perp),F_A(e_\alpha^\perp,e_\beta^\perp) \rangle \, w^\alpha \\ 
&+ \sum_{i=1}^{n-2} \sum_{\alpha=1}^2 \langle D_{A,e_i} D_{A,e_i} \phi,D_{A,e_\alpha^\perp} \phi 
\rangle \, w^\alpha \\ 
&= \varepsilon^2 \, \Big \langle d^{*_{g_0}} F_A - \frac{1}{2\varepsilon^2} \, (\phi \, 
\overline{D_A \phi} - \overline{\phi} \, D_A \phi),F_A(w,\cdot) \Big \rangle \\ 
&+ \Big \langle D_A^{*_{g_0}} D_A \phi - \frac{1}{2\varepsilon^2} \, (1 - |\phi|^2) \, 
\phi,D_{A,w} \phi \Big \rangle. 
\end{align*} 
We now take advantage of the identity 
\begin{align*} 
&\int_{NS_x} \sum_{\alpha,\beta=1}^2 \varepsilon^2 \, \langle F_A(e_i,e_\beta^\perp),F_A(e
_\alpha^\perp,e_\beta^\perp) \rangle \, w^\alpha \\ 
&+ \int_{NS_x} \sum_{\alpha=1}^2 \langle D_{A,e_i} \phi,D_{A,e_\alpha^\perp} \phi \rangle \, 
w^\alpha \\ 
&= -\frac{1}{2} \, \langle \nabla_i v,w \rangle \, \bigg ( 2 \int_{\mathbb{R}^2} \varepsilon^2 
\, |F_B|^2 + \int_{\mathbb{R}^2} |D_B \psi|^2 \bigg ) \\ 
&= -\frac{1}{2} \, \langle \nabla_i v,w \rangle \, \bigg ( \int_{\mathbb{R}^2} \varepsilon^2 
\, |F_B|^2 + \int_{\mathbb{R}^2} |D_B \psi|^2 + \int_{\mathbb{R}^2} \frac{1}{4\varepsilon^2} \, 
(1 - |\psi|^2)^2 \bigg ) \\ 
&= -\pi \, \langle \nabla_i v,w \rangle. 
\end{align*} 
Differentiating this identity, we obtain 
\begin{align*} 
&\int_{NS_x} \sum_{i=1}^{n-2} \sum_{\alpha,\beta=1}^2 \varepsilon^2 \, \nabla_{e_i} \langle F_A
(e_i,e_\beta^\perp),F_A(e_\alpha^\perp,e_\beta^\perp) \rangle \, w^\alpha \\ 
&+ \int_{NS_x} \sum_{i=1}^{n-2} \sum_{\alpha=1}^2 \nabla_{e_i} \langle D_{A,e_i} \phi,D_{A,e
_\alpha^\perp} \phi \rangle \, w^\alpha \\ 
&= \sum_{i=1}^{n-2} \nabla_{e_i} \int_{NS_x} \sum_{\alpha,\beta=1}^2 \varepsilon^2 \, \langle F
_A(e_i,e_\beta^\perp),F_A(e_\alpha^\perp,e_\beta^\perp) \rangle \, w^\alpha \\ 
&+ \sum_{i=1}^{n-2} \nabla_{e_i} \int_{NS_x} \sum_{\alpha=1}^2 \langle D_{A,e_i} \phi,D_{A,e
_\alpha^\perp} \phi \rangle \, w^\alpha \\ 
&- \int_{NS_x} \sum_{i=1}^{n-2} \sum_{\alpha,\beta=1}^2 \varepsilon^2 \, \langle F_A(e_i,e_\beta
^\perp),F_A(e_\alpha^\perp,e_\beta^\perp) \rangle \, \nabla_i w^\alpha \\ 
&- \int_{NS_x} \sum_{i=1}^{n-2} \sum_{\alpha=1}^2 \nabla_{e_i} \langle D_{A,e_i} \phi,D_{A,e
_\alpha^\perp} \phi \rangle \, \nabla_i w^\alpha \\ 
&= -\pi \, \langle \Delta v,w \rangle. 
\end{align*} 
Thus, we conclude that 
\begin{align*} 
&\int_{NS_x} \varepsilon^2 \, \Big \langle d^{*_{g_0}} F_A - \frac{1}{2\varepsilon^2} \, 
(\phi \, \overline{D_A \phi} - \overline{\phi} \, D_A \phi),F_A(w,\cdot) \Big \rangle \\ 
&+ \int_{NS_x} \Big \langle D_A^{*_{g_0}} D_A \phi - \frac{1}{2\varepsilon^2} \, (1 - 
|\phi|^2) \, \phi,D_{A,w} \phi \Big \rangle \\ 
&= -\pi \, \langle \Delta v,w \rangle. 
\end{align*} 
From this the assertion follows. \\

\begin{proposition}
The fibrewise projection of the error term to the obstruction space satisfies the estimate 
\begin{align*} 
\bigg \| 
&\Pi \Big ( d^* F_A - \frac{1}{2\varepsilon^2} \, (\phi \, \overline{D_A \phi} - \overline{\phi} 
\, D_A \phi),D_A^* D_A \phi - \frac{1}{2\varepsilon^2} \, (1 - |\phi|^2) \, \phi \Big ) \\ 
&- \Big ( \Delta v_\rho + \sum_{i,j=1}^{n-2} \sum_{\rho,\sigma=1}^2 h_{ij,\rho} \, h_{ij,\sigma} 
\, v_\sigma + \sum_{i=1}^{n-2} \sum_{\rho,\sigma=1}^2 R_{i\rho\sigma i} \, v_\sigma \Big ) \bigg 
\|_{\mathcal{C}^\gamma(S)} \leq C \, \varepsilon^2. 
\end{align*} 
\end{proposition}

\textbf{Proof.}
The Riemannian metric satisfies the asymptotic expansion 
\begin{align*} 
g(e_i,e_j) 
&= \delta_{ij} + 2 \sum_{\rho=1}^2 h_{ij,\rho} \, y_\rho \\ 
&+ \sum_{k=1}^{n-2} \sum_{\rho,\sigma=1}^2 h_{ik,\rho} \, h_{jk,\sigma} \, y_\rho \, y_\sigma \\ 
&- \sum_{\rho,\sigma=1}^2 R_{i\rho\sigma j} \, y_\rho \, y_\sigma + O(|y|^3) \\ 
g(e_i,e_\alpha^\perp) 
&= O(|y|^2) \\ 
g(e_\alpha^\perp,e_\beta^\perp) 
&= \delta_{\alpha\beta} - \frac{1}{3} \sum_{\rho,\sigma=1}^2 R_{\alpha\rho\sigma\beta} \, y_\rho 
\, y_\sigma + O(|y|^3). 
\end{align*} 
Using this asymptotic expansion, one can deduce Proposition 7.2 from Proposition 7.1. The 
details are left to the reader. \\

By Proposition 6.2, it suffices to find a normal vector field $v$ such that 
\[\Pi \Big ( d^* F_{\tilde{A}} - \frac{1}{2\varepsilon^2} \, (\tilde{\phi} \, \overline{D
_{\tilde{A}} \tilde{\phi}} - \overline{\tilde{\phi}} \, D_{\tilde{A}} \tilde{\phi}) + \frac
{1}{\varepsilon} \, du,D_{\tilde{A}}^* D_{\tilde{A}} \tilde{\phi} - \frac{1}{2\varepsilon^2} \, 
(1 - |\tilde{\phi}|^2) \, \tilde{\phi} - \frac{1}{\varepsilon} \, \tilde{\phi} \, u \Big ) = 0.
\] 
To this end, we need the following result. \\

\begin{proposition}
The pair $(\tilde{A},\tilde{\phi})$ satisfies 
\begin{align*} 
\bigg \| 
&\Pi \Big ( d^* F_{\tilde{A}} - \frac{1}{2\varepsilon^2} \, (\tilde{\phi} \, \overline{D_{\tilde
{A}} \tilde{\phi}} - \overline{\tilde{\phi}} \, D_{\tilde{A}} \tilde{\phi}) + \frac
{1}{\varepsilon} \, du,D_{\tilde{A}}^* D_{\tilde{A}} \tilde{\phi} - \frac{1}{2\varepsilon^2} \, 
(1 - |\tilde{\phi}|^2) \, \tilde{\phi} - \frac{1}{\varepsilon} \, \tilde{\phi} \, u \Big ) \\ 
&- \Big ( \Delta v_\rho + \sum_{i,j=1}^{n-2} \sum_{\rho,\sigma=1}^2 h_{ij,\rho} \, h_{ij,\sigma} 
\, v_\sigma + \sum_{i=1}^{n-2} \sum_{\rho,\sigma=1}^2 R_{i\rho\sigma i} \, v_\sigma \Big ) \bigg 
\|_{\mathcal{C}^\gamma(S)} \leq C \, \varepsilon^{\frac{5}{4}}. 
\end{align*} 
\end{proposition}

\textbf{Proof.} 
Using the estimate 
\[\|(a,f)\|_{\mathcal{C}_{\mu,\varepsilon}^{2,\gamma}(M)} \leq C \, \varepsilon^2,\] 
we obtain 
\[\|\Pi \, Q(a,f)\|_{\mathcal{C}^\gamma(S)} \leq C \, \varepsilon^{1-\gamma} \, \|Q(a,f)\|
_{\mathcal{C}_{\mu,\varepsilon}^\gamma(M)} \leq C \, \varepsilon^{3-\gamma}.\] 
Moreover, we have 
\[\|\Pi \, \mathbb{L}_\varepsilon (a,f)\|_{\mathcal{C}^\gamma(S)} \leq C \, \varepsilon^{\frac
{5}{4}}.\] 
The assertion follows now from Proposition 7.2. \\

\begin{proposition}
If $S$ is non-degenerate, then we can choose the normal vector field $v$ such that 
\[\Pi \Big ( d^* F_{\tilde{A}} - \frac{1}{2\varepsilon^2} \, (\tilde{\phi} \, \overline{D
_{\tilde{A}} \tilde{\phi}} - \overline{\tilde{\phi}} \, D_{\tilde{A}} \tilde{\phi}) + \frac
{1}{\varepsilon} \, du,D_{\tilde{A}}^* D_{\tilde{A}} \tilde{\phi} - \frac{1}{2\varepsilon^2} \, 
(1 - |\tilde{\phi}|^2) \, \tilde{\phi} - \frac{1}{\varepsilon} \, \tilde{\phi} \, u \Big ) = 0.
\] 
Therefore, the pair $(\tilde{A},\tilde{\phi})$ corresponding to this normal vector field $v$ is 
a solution of the Ginzburg-Landau equations.
\end{proposition}

\textbf{Proof.}
Let $J$ be the Jacobi-operator of the submanifold $S$. According to Proposition 7.3, we may 
write 
\begin{align*} 
&\Pi \Big ( d^* F_{\tilde{A}} - \frac{1}{2\varepsilon^2} \, (\tilde{\phi} \, \overline{D_{\tilde
{A}} \tilde{\phi}} - \overline{\tilde{\phi}} \, D_{\tilde{A}} \tilde{\phi}) + \frac
{1}{\varepsilon} \, du,D_{\tilde{A}}^* D_{\tilde{A}} \tilde{\phi} - \frac{1}{2\varepsilon^2} \, 
(1 - |\tilde{\phi}|^2) \, \tilde{\phi} - \frac{1}{\varepsilon} \, \tilde{\phi} \, u \Big ) \\ 
&= Jv + R(v), 
\end{align*} 
where $\|R(v)\|_{\mathcal{C}^\gamma(S)} \leq C \, \varepsilon^{\frac{5}{4}}$ for $\|v\|
_{\mathcal{C}^{2,\gamma}(S)} \leq \varepsilon$. Hence, the mapping $-J^{-1} \, R$ maps a ball of 
radius $\varepsilon$ in the Banach space $\mathcal{C}^{2,\gamma}(S)$ into a ball of radius $C \, 
\varepsilon^{\frac{5}{4}}$ in $\mathcal{C}^{2,\gamma}(S)$. Unfortunately, Schauder's fixed point 
theorem cannot be applied, since the mapping $-J^{-1} \, R$ need not be compact. To overcome 
this problem, we use an idea of F. Pacard and M. Ritor\'e \cite{PR1}. Using an appropriate 
sequence of smoothing operators, we may approximate the mapping $-J^{-1} \, R$ by a sequence of 
compact mappings. Each of these mappings has a fixed point in $\mathcal{C}^{2,\gamma}(S)$ by 
Schauder's fixed point theorem. Taking limits, we obtain a fixed point of the original mapping 
$-J^{-1} \, R$ in the Banach space $\mathcal{C}^{2,\frac{\gamma}{2}}(S)$. Hence, there exists a 
normal vector field $v \in \mathcal{C}^{2,\frac{\gamma}{2}}(S)$ such that $Jv + R(v) = 0$. 
Hence, the pair $(\tilde{A},\tilde{\phi})$ corresponding to that choice of the vector field $v$ 
satisfies 
\[\Pi \Big ( d^* F_{\tilde{A}} - \frac{1}{2\varepsilon^2} \, (\tilde{\phi} \, \overline{D
_{\tilde{A}} \tilde{\phi}} - \overline{\tilde{\phi}} \, D_{\tilde{A}} \tilde{\phi}) + \frac
{1}{\varepsilon} \, du,D_{\tilde{A}}^* D_{\tilde{A}} \tilde{\phi} - \frac{1}{2\varepsilon^2} \, 
(1 - |\tilde{\phi}|^2) \, \tilde{\phi} - \frac{1}{\varepsilon} \, \tilde{\phi} \, u \Big ) = 0.
\] 
This concludes the proof. \\

\end{document}